\documentclass[12pt,a4paper,draft]{article}%
\usepackage{graphicx}
\usepackage{amsmath}
\usepackage{amsfonts}
\usepackage{amssymb}%
\providecommand{\U}[1]{\protect\rule{.1in}{.1in}}
\providecommand{\U}[1]{\protect\rule{.1in}{.1in}}
\newtheorem{theorem}{Theorem}

\newtheorem{corollary}[theorem]{Corollary}

\newtheorem{definition}[theorem]{Definition}

\newtheorem{lemma}[theorem]{Lemma}

\newtheorem{proposition}[theorem]{Proposition}
\newtheorem{remark}[theorem]{Remark}

\newenvironment{proof}[1][Proof]{\noindent\textbf{#1.} }{\ \rule{0.5em}{0.5em}}
\begin{document}

\title{Multidimensional Chebyshev spaces, hierarchy of infinite-dimensional spaces
and Kolmogorov-Gelfand widths}
\author{O. Kounchev\\Institute of Mathematics and Informatics, \\Bulgarian Academy of Sciences\\and\\IZKS, University of Bonn}
\maketitle

\begin{abstract}
Recently the theory of widths of Kolmogorov (especially of Gelfand widths) has
received a great deal of interest due to its close relationship with the newly
born area of Compressed Sensing. It has been realized that widths reflect
properly the sparsity of the data in Signal Processing. However fundamental
problems of the theory of widths in multidimensional Theory of Functions
remain untouched, and their progress will have a major impact over analogous
problems in the theory of multidimensional Signal Analysis. The present paper
has three major contributions:\ 

1. We solve the longstanding problem of finding multidimensional
generalization of the Chebyshev systems: we introduce Multidimensional
Chebyshev spaces, based on solutions of higher order elliptic equation, as a
generalization of the one-dimensional Chebyshev systems, more precisely of the
\emph{ECT--}systems.

2. Based on that we introduce a new hierarchy of infinite-dimensional spaces
for functions defined in multidimensional domains; we define corresponding
generalization of Kolmogorov's widths.

3. We generalize the original results of Kolmogorov by computing the widths
for special "ellipsoidal" sets of functions defined in multidimensional domains.

\textbf{MSC2010: 35J, 46L, 41A, 42B.}

\textbf{Keywords:} Kolmogorov widths, Gelfand widths, Compressed Sensing,
Chebyshev systems, Approximation by solutions of PDEs.

\end{abstract}

\section{Introduction}

It is a notorious fact that the polynomials of several variables fail to enjoy
the nice interpolation and approximation properties of the one-dimensional
polynomials, and this is particularly visible in such fundamental areas of
Mathematical Analysis as the Moment problem, Interpolation, Approximation,
etc. One alternative approach is to use solutions of elliptic equations, in
particular polyharmonic functions, which has led to new amazing Ansatzes in
multidimensional Approximation and Interpolation \cite{hayman},
\cite{goldstein}, \cite{fuglede}, \cite{kounchev1992TAMS},
\cite{haymanKorenblum}, in Spline and Wavelet Theory \cite{madych},
\cite{okbook}, and recently in the Moment Problem \cite{kounchevrenderMoment}.
This approach has been given the name \emph{Polyharmonic Paradigm}
\cite{okbook}, as an approach to Multidimensional Mathematical Analysis, which
is opposite to the usual concept which is based on algebraic and trigonometric
polynomials of several variables. However the effectiveness of the
Polyharmonic Paradigm has remained unexplained for a long time.

One of the \emph{main objectives} of the present research is to find a new
point of view on the longstanding problem of finding multidimensional
generalization of Chebyshev systems: In Section \ref{SIntroHierarchy} we show
that the solution spaces of a wide class of elliptic PDEs (defined further as
\emph{Multidimensional Chebyshev spaces}) are a natural generalization of the
one-variable polynomials as well as of the one-dimensional Chebyshev
systems.\footnote{Recall that the one-dimensional Chebyshev systems appeared
as a generalization of the one-variable algebraic polynomials in the works of
A. Markov in the context of the classical Moment problem. They were further
developed and applied to the generalized Moment problem and Approximation
theory by A. Haar (the Haar spaces), S. Bernstein, M. Krein, S. Karlin, and
others, cf. \cite{kreinnudelman}, \cite{karlinstudden}, \cite{mccullough}. }
In particular, this shows that  the polyharmonic functions which are solutions
to the polyharmonic operators are a generalization of the one-variable
polynomials. Let us recall that there  has been a long search for proper
multidimensional generalization of Chebyshev systems. The standard
generalization by means of zero set property fails to produce a non-trivial
multidimensional system which is the content of the theorem of Mairhuber, see
the thorough discussion in \cite{kreinnudelman} (chapter $2,$ section $1.1$).
Our generalization provided by the Multidimensional Chebyshev spaces comes
from a completely different perspective, by generalizing the boundary value
properties of the one-dimensional polynomials.\footnote{In general, zero set
properties and intersections are not a reliable reference point for
multidimensional Analysis. In particular, let us recall that polyharmonic (and
even harmonic) functions do not have simple zero sets, although they are
solutions to nice BVPs as (\ref{BVP1})-(\ref{BVP2}), cf.
\cite{haymanKorenblum}.}

On the other hand the \emph{cornerstone} of the present paper is an amazing
though simple characterization of the $N-$dimensional subspaces $X_{N}$ of
$C^{N-1}\left(  I\right)  $ (here $I$ is an interval in $\mathbb{R}$) via
Chebyshev systems. This characterization says that the typical subspace
$X_{N}$ is a finite-piecewise Chebyshev space, or to be more correct,
finite-piecewise $\emph{ECT-}$system. This discovery causes an immediate chain
reaction: by analogy, for a domain $D\subset\mathbb{R}^{n},$ we use the
newly-invented Multidimensional Chebyshev spaces to define in $L_{2}\left(
D\right)  $ a multidimensional generalization $\mathcal{X}_{N}$ of the spaces
$X_{N},$ which we call "spaces of Harmonic Dimension $N$". These spaces
$\mathcal{X}_{N}$ represent a \emph{new hierarchy of infinite-dimensional
spaces. }Hence, the big surprise of the present research is the
reconsideration of the simplistic understanding that the natural
generalization of $X_{N}$ is provided just by the finite-dimensional subspaces
of $L_{2}\left(  D\right)  .$ We realize that the finite-dimensional subspaces
in $C^{N}\left(  D\right)  ,$ for domains $D\subset\mathbb{R}^{n}$ for
$n\geq2,$ \textbf{do not} serve the same job as the finite-dimensional
subspaces in $C^{N}\left(  D\right)  $ for intervals $D\subset\mathbb{R}^{1},$
and one has to replace them by a lot more sophisticated objects, namely by the
spaces having Harmonic Dimension $N.$

Respectively, \emph{the focus} of the present research is, by means of the
spaces $\mathcal{X}_{N}$ to define a multidimensional generalization of the
Kolmogorov-Gelfand $N-$widths, which we call "Harmonic $N-$widths". After that
we compute the Harmonic $N-$widths for "cylindrical ellipsoids" in
$L_{2}\left(  D\right)  ,$ by generalizing the original results of Kolmogorov.

Another important motivation for the present research is the recent interest
to the theory of widths (especially to Gelfand widths) coming from the
applications in an area of Signal Analysis, called Compressed Sensing (CS). In
a certain sense the central idea of CS is rooted in the theory of widths, cf.
e.g. \cite{donoho}, \cite{candes}, \cite{devore}, \cite{pinkus2011}. However,
apparently this strategy works smoothly only in the case of representation of
one-dimensional signals, while an adequate approach to multivariate signals is
missing -- one reason may be found by analogy in the fact that the theory of
Kolmogorov-Gelfand widths fits properly only for one-dimensional function
spaces (as pointed out below, e.g. in formula (\ref{dNisInfinity})). Recently,
a new multivariate Wavelet Analysis was developed based on solutions of
elliptic partial differential equations (\cite{okbook}), in particular
"polyharmonic subdivision wavelets" were introduced (cf.
\cite{dynkounchevlevinrender}, \cite{kounchevKalaglarsky}). In order to apply
CS ideas to these wavelets it would require essential generalization of the
theory of widths for infinite-dimensional spaces, and it is expected that the
present research is a step in the right direction.

In his seminal paper \cite{kolmogorov1936} Kolmogorov has introduced the
theory of widths and has applied it ingeniously to the following set of
functions defined in the compact interval:
\begin{equation}
K_{p}:=\left\{  f\in AC^{p-1}\left(  \left[  a,b\right]  \right)  :%
{\displaystyle\int_{0}^{1}}
\left\vert f^{\left(  p\right)  }\left(  t\right)  \right\vert ^{2}%
dt\leq1\right\}  . \label{Kp}%
\end{equation}
In the present paper we study a natural multivariate generalization of the set
$K_{p}$ which in a domain $B\subset\mathbb{R}^{n}$ is given by
\begin{equation}
K_{p}^{\ast}:=\left\{  u\in H^{2p}\left(  B\right)  :%
{\displaystyle\int_{B}}
\left\vert \Delta^{p}u\left(  x\right)  \right\vert ^{2}dx\leq1\right\}  ,
\label{Kpstar}%
\end{equation}
where $\Delta^{p}$ is the $p-$th iterate of the Laplace operator $\Delta=%
{\displaystyle\sum_{j=1}^{n}}
\partial^{2}/\partial x_{j}^{2};$ we consider more general sets $K_{p}^{\ast}$
given in (\ref{KpstarGeneral}) below.

Let us summarize the major contributions of the present paper:

\begin{enumerate}
\item For every integer $N\geq0$ we define the \emph{Multidimensional
Chebyshev spaces} of order $N$ as spaces of solutions of a special class of
elliptic PDEs of order $2N.$ They generalize the classical one-dimensional
Extended Complete Chebyshev systems (\emph{ECT--}systems).

\item For every integer $N\geq0$ we generalize the $N-$dimensional subspaces
$X_{N}$ in $C^{N-1}\left(  I\right)  $ (for intervals $I\subset\mathbb{R}$) to
a multidimensional setting. The generalization $\mathcal{X}_{N}$ is a
piecewise Multidimensional Chebyshev \allowbreak space of order $N,$ and is
said to have "Harmonic Dimension $N$". This represents a new hierarchy of
infinite-dimensional spaces of functions defined in domains in $\mathbb{R}%
^{n}.$

\item For every integer $N\geq0$ we define \textbf{Harmonic Widths} which
generalize the Kolmogorov widths, whereby we use as approximants the spaces
$\mathcal{X}_{N}$ instead of finite-dimensional spaces $X_{N}$ used by the
Kolmogorov widths. We generalize the one-dimensional Kolomogorov's results in
the theory of widths.
\end{enumerate}

The \emph{crux} of the new notion of hierarchy of infinite-dimensional spaces
is the following: Let the domain $D\subset\mathbb{R}^{n}$ be compact with
sufficiently smooth boundary $\partial D.$ Then the $N-$dimensional subspaces
in $C^{\infty}\left(  I\right)  $ will be generalized by spaces of
\emph{solutions of elliptic equations} (and by more general spaces introduced
in Definition \ref{Dhdimension} below):
\begin{equation}
\mathcal{X}_{N}=\left\{  u:P_{2N}u\left(  x\right)  =0,\quad\text{for }x\in
D\right\}  \subset L_{2}\left(  D\right)  ; \label{XN}%
\end{equation}
here $P_{2N}$ is an elliptic operator of order $2N$ in the domain $D.$
Respectively, the simplest version of our generalization of Kolmogorov's
theorem about widths finds the extremizer of the following problem
\[
\inf_{\mathcal{X}_{N}}\operatorname*{dist}\left(  \mathcal{X}_{N},K_{p}^{\ast
}\right)  ,
\]
where $K_{p}^{\ast}$ is the set defined in (\ref{Kpstar}) and $\mathcal{X}%
_{N}$ is defined in (\ref{XN}) by an elliptic operator $P_{2N}$ of order $2N;$
for the complete formulation see Theorem \ref{TKolmogorovMultivariate} below.

What is the reason to take namely solutions of elliptic equations in the
multidimensional case is explained in the following section.

\textbf{Acknowledgements:} The author acknowledges the support of the
\allowbreak\ Alexander von Humboldt Foundation, and of Project DO-02-275 with
Bulgarian NSF. The author thanks the following Professors: Matthias Lesch for
the interesting discussion about hierarchies of infinite-dimensional linear
spaces, Hermann Render about advice on multivariate polynomial division, and
Peter Popivanov, Nikolay Kutev and Georgi Boyadzhiev about advice on Elliptic BVP.

\section{Multidimensional Chebyshev spaces and a hierarchy of
infinite-dimensional spaces \label{SIntroHierarchy}}

Let us give a heuristic outline of \emph{the main idea} of this new hierarchy
of spaces, by explaining how it appears as a natural generalization of the
finite-dimensional subspaces of $C^{N-1}\left[  a,b\right]  $ in a compact
interval $\left[  a,b\right]  $ in $\mathbb{R}.$ First of all, we will show
that there exists an amazing relation between the finite-dimensional subspaces
of $C^{N-1}\left[  a,b\right]  $ and the theory of Chebyshev systems.

Let us construct some special type of $N-$dimensional subspaces of
\allowbreak\ $C^{N-1}\left(  J\right)  $ where the interval $J=\left(
c,d\right)  \subset\mathbb{R}.$ Let the functions $\rho_{j},$ $j=1,2,...,N$
satisfy $\rho_{j}>0$ on $J,$ and $\rho_{j}\in C^{N+1-j}\left(  J\right)  .$ We
assume that the functions $\rho_{j}$ satisfy the necessary integrability so
that we may define the following functions:
\begin{align}
v_{1}\left(  t\right)   &  =\rho_{1}\left(  t\right) \label{explicitECT1}\\
v_{2}\left(  t\right)   &  =\rho_{1}\left(  t\right)
{\displaystyle\int_{c}^{t}}
\rho_{2}\left(  t_{2}\right)  dt_{2}\\
&  \cdot\cdot\cdot\nonumber\\
v_{N}\left(  t\right)   &  =\rho_{1}\left(  t\right)
{\displaystyle\int_{c}^{t}}
\rho_{2}\left(  t_{2}\right)
{\displaystyle\int_{c}^{t_{2}}}
\rho_{3}\left(  t_{3}\right)  \cdot\cdot\cdot%
{\displaystyle\int_{c}^{t_{N-1}}}
\rho_{N}\left(  t_{N}\right)  dt_{2}dt_{3}\cdot\cdot\cdot dt_{N}.
\label{explicitECT4}%
\end{align}

Let us recall some classical results about the space $X_{N}%
=\operatorname*{span}\left\{  v_{j}\right\}  _{j=1}^{N}.$ For every
$k=1,2,...,N$ the consecutive Wronskians for the system of functions $\left\{
v_{j}\right\}  _{j=1}^{N}$ may be computed explicitly and are given by%
\begin{equation}
W_{k}:=W\left(  v_{1},v_{2},...,v_{k}\right)  =\rho_{1}^{k}\rho_{2}^{k-1}%
\cdot\cdot\cdot\rho_{k}, \label{Wro}%
\end{equation}
and vice versa:
\begin{align}
\rho_{1}  &  =W_{1}=v_{1},\quad\rho_{2}=W_{2}/W_{1}^{2}\label{roW1}\\
\rho_{k}  &  =W_{k}W_{k-2}/W_{k-1}^{2}\qquad\text{for }k\geq3, \label{roW2}%
\end{align}
(cf. \cite{polya}, or \cite{karlinstudden}, chapter $11,$ formulas $1.12$ and
$1.13$). From representation formula (\ref{Wro}) directly follows that for all
$k=1,2,...,N$ the Wronskians satisfy
\begin{equation}
W\left(  v_{1},v_{2},...,v_{k}\right)  >0\qquad\text{on }J.
\label{WronskiSign}%
\end{equation}

Let us define on $J$ the ordinary differential operator
\begin{equation}
L_{N}\left(  t;\frac{d}{dt}\right)  =\frac{d}{dt}\frac{1}{\rho_{N}\left(
t\right)  }\cdot\cdot\cdot\frac{d}{dt}\frac{1}{\rho_{2}\left(  t\right)
}\frac{d}{dt}\frac{1}{\rho_{1}\left(  t\right)  }, \label{LNoperatorForm}%
\end{equation}
Then, from formulas (\ref{explicitECT1})--(\ref{explicitECT4}) directly
follows that all $v_{j}$'s, hence all elements of the space
$\operatorname*{span}\left\{  v_{j}\right\}  _{j=1}^{N},$ satisfy the ODE
\begin{equation}
L_{N}u\left(  t\right)  =0,\qquad\text{for }t\in J. \label{LN+1}%
\end{equation}
Obviously, the operator $L_{N}$ has a non-negative leading coefficient and is
in this sense one-dimensional "elliptic".

We have the following classical result (cf. \cite{karlinstudden}, chapter
$11,$ Theorem $1.2$).

\begin{proposition}
The space
\begin{equation}
X_{N}=\operatorname*{span}\left\{  v_{j}\right\}  _{j=1}^{N} \label{XN1dim}%
\end{equation}
is an $N-$dimensional subspace of $C^{N-1}\left(  J\right)  .$
\end{proposition}

We recall the following definition (cf. \cite{karlinstudden}, chapter $11$).

\begin{definition}
A space $X_{N}\subset C^{N-1}\left(  J\right)  $ is called \textbf{ECT--space}
(or \textbf{Extended Complete Chebyshev space}) if it has a basis $\left\{
v_{j}\right\}  _{j=1}^{N}$ satisfying the Wronskian condition
(\ref{WronskiSign}).
\end{definition}

Our terminology above differs slightly from the accepted terminology:\ we
consider Chebyshev systems on open intervals $J$ and instead of "$\emph{ECT-}%
$system" we say "$\emph{ECT-}$space".

The following result characterizes the typical ("general position")
$N-$dimensional subspaces of \allowbreak$C^{N-1}\left[  a,b\right]  $ by means
of the $\emph{ECT-}$spaces.

\begin{theorem}
\label{PstructureMarkovSystem} "Almost all" $N-$dimensional subspaces of
$C^{N-1}\left[  a,b\right]  $ are finite-piecewise \emph{ECT--}spaces in the
following sense: If $X_{N}$ is an $N-$dimensional subspace of $C^{N-1}\left[
a,b\right]  $ then there exists a sequence of $N-$dimensional subspaces
$X_{N}^{m}$ of $C^{N-1}\left[  a,b\right]  ,$ $m\geq1,$ satisfying:

1. For every space $X_{N}^{m}$ there exists a finite subdivision $a\leq
t_{0}<t_{1}<\cdot\cdot\cdot<t_{p_{m}}=b$ of the interval $\left[  a,b\right]
$ such that the restriction of $X_{N}^{m}$ to every subinterval $\left(
t_{j},t_{j+1}\right)  $ is an \emph{ECT--}space.

2. The following limit holds in the metric of $C^{N-1}\left[  a,b\right]  ,$
\[
U\left(  X_{N}^{m}\right)  \longrightarrow U\left(  X_{N}\right)
\qquad\text{for }m\longrightarrow\infty,
\]
where $U$ denotes the unit ball in the corresponding space.
\end{theorem}

The proof and an example are provided in the Appendix in Section
\ref{Sappendix1} below.

Hence, from Theorem \ref{PstructureMarkovSystem} and formula (\ref{LN+1}) we
see that a typical $N-$dimen\-sio\-nal subspace of $C^{N-1}\left[  a,b\right]
$ is a piecewise $\emph{ECT-}$space, i.e. piecewise solution space of a family
of ordinary differential operators $L_{N}$ (with coefficients depending on the
intervals $\left(  t_{j},t_{j+1}\right)  $ ). In our multidimensional
generalization we will generalize these typical $N-$dimensional subspaces of
$C^{N-1}\left[  a,b\right]  .$ Indeed, Theorem \ref{PstructureMarkovSystem}
already suggests an \textbf{Ansatz for Multidimensional Chebyshev spaces }
which generalize the one-dimensional \emph{ECT--}spaces: First, we generalize
the one-dimensional operators $L_{N}$ in (\ref{LNoperatorForm}) by considering
in a domain $D\subset\mathbb{R}^{n}$ \emph{elliptic partial differential
operators}\textbf{ }of the form
\[
P_{2N}\left(  x,D_{x}\right)  =Q^{\left(  1\right)  }Q^{\left(  2\right)
}\cdot\cdot\cdot Q^{\left(  N\right)  }%
\]
where $Q^{\left(  j\right)  }$ are \emph{elliptic operators of second order}
in $D.$ Then the corresponding solution space $\mathcal{X}_{N}$ defined by
\begin{equation}
\mathcal{X}_{N}:=\left\{  u\in C^{\infty}\left(  D\right)  :P_{2N}u\left(
x\right)  =0\quad\text{in }D\right\}  , \label{S2M}%
\end{equation}
is our generalization of the one-dimensional \emph{ECT--}space, and we call it
\textbf{Multidimensional Chebyshev space.}

However, we need to impose some more conditions on the elliptic operators
$Q^{\left(  j\right)  }:$ We would like that the elements of the space
$\mathcal{X}_{N}$ generalize the \emph{interpolation properties} of the
one-dimensional \emph{ECT--}spaces $X_{N},$ and this would require more
conditions to be imposed on the operators $Q^{\left(  j\right)  }.$ Complete
analogy between the one-dimensional and the multidimensional case is achieved
only for the one-dimensional \emph{ECT--}spaces of\emph{ even dimension} which
satisfy the following interpolation property (\cite{karlinstudden}, chapter
$11$):

\begin{proposition}
\label{PDirichletBVP} Let the \emph{ECT--}space $X_{2M}\subset C^{2M-1}\left(
I\right)  $ be given. Then for every subinterval $I_{1}=\left[  a_{1}%
,b_{1}\right]  \subset I,$ and for arbitrary constants $\left\{  c_{k}%
,d_{k}\right\}  _{k=0}^{M-1},$ the (Dirichlet) boundary value problem
\begin{align}
u^{\left(  k\right)  }\left(  a_{1}\right)   &  =c_{k}\qquad\text{for
}k=0,1,...,M-1\label{Dirichlet1}\\
u^{\left(  k\right)  }\left(  b_{1}\right)   &  =d_{k}\qquad\text{for
}k=0,1,...,M-1 \label{Dirichlet2}%
\end{align}
has a solution $u\in X_{2M}.$
\end{proposition}

The interpolation property of Proposition \ref{PDirichletBVP} reminds us
immediately of the solvability of the Elliptic BVP. We specify below the
well-known Dirichlet BVP for the operator $P_{2N}$ considered on subdomains
$D_{1}\subset D$:
\begin{align}
P_{2N}u\left(  x\right)   &  =0\quad\text{for }x\in D_{1}%
\label{DirichletMultivariate}\\
\left(  \frac{\partial}{\partial n}\right)  ^{k}u\left(  y\right)   &
=c_{k}\left(  y\right)  \qquad\text{for }y\in\partial D_{1},\quad\text{for
}k=0,1,...,N-1. \label{DirichletMultivariate2}%
\end{align}
Thus the solvability of the one-dimensional problem (\ref{Dirichlet1}%
)-(\ref{Dirichlet2}) in the space $X_{2M}$ may be considered as a special case
of the multidimensional theory for \textbf{Elliptic Boundary Value Problems}
(\textbf{BVP}).

Let us remind that the Dirichlet BVP in (\ref{DirichletMultivariate}%
)-(\ref{DirichletMultivariate2}) is well-known to be solvable for data
$\left\{  c_{k}\left(  y\right)  \right\}  _{k=0}^{N-1}$ from a proper Sobolev
or H\"{o}lder space on the boundary $\partial D_{1}.$ An important point is
that for a large class of elliptic operators $P_{2N}$ every solution of
(\ref{DirichletMultivariate})-(\ref{DirichletMultivariate2}) may be
approximated by solutions in the whole domain $D,$ i.e. by elements of
$\mathcal{X}_{N}.$ This may be considered as a substitute of the interpolation
property (\ref{Dirichlet1})-(\ref{Dirichlet2}) in the one-dimensional case.
This is also the explanation for the judicious choice of the special class of
operators $P_{2N}$ in Definition \ref{Dhdimension} below, as we mimic the
operators in (\ref{LNoperatorForm}) by satisfying some natural interpolation properties.

Making analogy with the one-dimensional case (\ref{Dirichlet1}%
)-(\ref{Dirichlet2}), we may say that here the space $\mathcal{X}_{N}$ defined
in (\ref{S2M}) is "parametrized" by the boundary conditions $B_{j}u,$ however
the "parameters" $\left\{  c_{k}\left(  y\right)  \right\}  _{k=0}^{M-1}$ run
a function space. Hence, the spaces $\mathcal{X}_{N}$ may be considered as a
natural generalization of the one-dimensional \emph{ECT--}spaces and for that
reason we call them \textbf{Multidimensional Chebyshev spaces}.

After having defined the Multidimensional Chebyshev spaces, the next step will
be to introduce the promised multidimensional generalization of the "typical"
$N-$dimensional subspaces of $C^{N-1}\left(  I\right)  .$ We will define them
in Definition \ref{Dhdimension} below as subspaces $X_{N}$ of functions in
$L_{2}\left(  D\right)  $ which are piecewise solutions of (regular) elliptic
differential operators $P_{2N}$ of order $2N.$ We will say that $X_{N}$ has
"Harmonic Dimension $N$" and we will write
\[
\operatorname{hdim}\left(  X_{N}\right)  =N
\]
Kolmogorov's notion of $N-$width (and in a similar way Gelfand's width) is
naturally generalized for symmetric sets by the notion of "Harmonic
$\allowbreak N-$width" defined by putting
\[
\operatorname*{hd}\nolimits_{N}\left(  S\right)  :=\inf_{\operatorname{hdim}%
\left(  X_{N}\right)  =N}\operatorname*{dist}\left(  X_{N},S\right)  ,
\]
see Definition \ref{Dwidth} below. The \textbf{main result} of the present
paper is the computation of
\[
\operatorname*{hd}\nolimits_{N}\left(  K_{p}^{\ast}\right)  \qquad\text{for
}N\leq p,
\]
where $K_{p}^{\ast}$ is defined in (\ref{Kpstar}) and more generally in
(\ref{KpstarGeneral}).

\section{Plan of the paper}

To facilitate the reader, in section \ref{SKolmogorovsResults} we provide a
short summary of the original Kolmogorov's results. For the same reason, in
section \ref{SreminderBVP} we provide a short reminder on Elliptic BVP. In
section \ref{SprincipalAxes} we prove the representation of the "cylindrical
ellipsoid" set $K_{p}^{\ast}$ in principal axes which generalizes the
one-dimensional representation of Kolmogorov, cf. Theorem
\ref{TprincipalAxesMultivariate} below. In section \ref{Shierarchy} we
introduce the notion of \emph{Harmonic Dimension}, and the \emph{First Kind}
spaces of Harmonic Dimension $N.$ Based on it we define \emph{Harmonic Widths}
which generalize Kolmogorov's widths. In section \ref{Swidths}, in Theorem
\ref{TKolmogorovMultivariate} we prove a genuine analog to Kolmogorov's
theorem about widths. It says that among all spaces $X_{N}$ having Harmonic
Dimension $N,$ some special space $\widetilde{X}_{N}$ provides the best
approximation to the set $K_{p}^{\ast}$ in problem
\[
\inf_{X_{N}}\operatorname*{dist}\left(  X_{N},K_{p}^{\ast}\right)  ,
\]
and this space $\widetilde{X}_{N}$ is identified by the principal axes
representation provided by Theorem \ref{TprincipalAxesMultivariate}. In
section \ref{SsecondKind} we introduce \emph{Second Kind} spaces of Harmonic
Dimension $N$ and formulate a further generalization of Theorem
\ref{TKolmogorovMultivariate}. Apparently, the First and Second Kind spaces
having Harmonic Dimension $N$ provide the maximal generalization in the
present framework.

A special case of the present results is available in \cite{kounchevSozopol},
and might be instructive for the reader to start with.

A final remark to our generalization is in order. In our consideration we will
not strive to achieve a maximal generality. As it is clear, especially in the
applications to the theory of widths even in the one-dimensional case we may
consider not all $N-$dimensional subspaces but "almost all" $N-$dimensional
subspaces of $C^{\infty}\left(  D\right)  $ in some sense, or a class of
$N-$dimensional subspaces which are dense (in a proper topology) in the set of
all other $N-$dimensional subspaces. This "genericity" point of view is
essential in our multivariate generalization since it will allow us to avoid
burdensome proofs necessary in the case of the bigger generality of the
construction. For the same reason we will not consider elliptic
pseudo-differential operators although almost all results have a
generalization for such setting.

\section{Kolmogorov's results - a reminder \label{SKolmogorovsResults}}

In order to make our multivariate generalization transparent we will recall
the original results of Kolmogorov provided in his seminal paper
\cite{kolmogorov1936}. Kolmogorov has considered the set $K_{p}$ defined in
(\ref{Kp}). He proved that this is an \textbf{ellipsoid} by constructing
explicitly its principal axes. Namely, he considered the eigenvalue problem
\begin{align}
\left(  -1\right)  ^{p}u^{\left(  2p\right)  }\left(  t\right)   &  =\lambda
u\left(  t\right)  \qquad\qquad\qquad\text{for }t\in\left(  0,1\right)
\label{eigen1}\\
u^{\left(  p+j\right)  }\left(  0\right)   &  =u^{\left(  p+j\right)  }\left(
1\right)  =0\qquad\text{for }j=0,1,...,p-1. \label{eigen2}%
\end{align}
Kolmogorov used the following properties of problem (\ref{eigen1}%
)-(\ref{eigen2}) (cf. \cite{lorentz}, Chapter $9.6,$ Theorem $9,$ p. $146,$ or
\cite{naimark}, \cite{pinkus}):

\begin{proposition}
\label{PBVP}Problem (\ref{eigen1})-(\ref{eigen2}) has a countable set of
non-negative real eigenvalues with finite multiplicity. If we denote them by
$\lambda_{j}$ in a monotone order, they satisfy $\lambda_{j}\longrightarrow
\infty$ for $j\longrightarrow\infty.$ They satisfy the following asymptotic
$\lambda_{j}=\pi^{2p}j^{2p}\left(  1+O\left(  j^{-1}\right)  \right)  .$ The
corresponding orthonormalized eigenfunctions $\left\{  \psi_{j}\right\}
_{j=1}^{\infty}$ form a complete orthonormal system in $L_{2}\left(  \left[
0,1\right]  \right)  .$ The eigenvalue $\lambda=0$ has multiplicity $p$ and
the corresponding eigenfunctions $\left\{  \psi_{j}\right\}  _{j=1}^{p}$ are
the basis for the solutions to equation $u^{\left(  p\right)  }\left(
t\right)  =0$ in the interval $\left(  0,1\right)  .$
\end{proposition}

Further, Kolmogorov provided a description of the axes of the "cylindrical
ellipsoid" $K_{p}$, from which an approximation theorem of \emph{Jackson type
}easily follows (cf. \cite{lorentz}, chapter $4$ and chapter $5$).

\begin{proposition}
\label{PKolmogorovJackson} Let $f\in L_{2}\left(  \left[  a,b\right]  \right)
$ have the $L_{2}-$expansion
\[
f\left(  t\right)  =%
{\displaystyle\sum_{j=1}^{\infty}}
f_{j}\psi_{j}\left(  t\right)  .
\]
Then $f\in K_{p}$ if and only if
\[%
{\displaystyle\sum_{j=1}^{\infty}}
f_{j}^{2}\lambda_{j}\leq1.
\]
For $N\geq p+1$ and every $f\in K_{p}$ holds the following estimate
(\emph{Jackson type} approximation):
\begin{equation}
\left\Vert f-%
{\displaystyle\sum_{j=1}^{N}}
f_{j}\psi_{j}\left(  t\right)  \right\Vert _{L_{2}}\leq\frac{1}{\sqrt
{\lambda_{N+1}}}=O\left(  \frac{1}{\left(  N+1\right)  ^{p}}\right)  .
\label{Jackson1dim}%
\end{equation}

\end{proposition}

However, Kolmogorov didn't stop at this point but asked further, whether the
linear space $\widetilde{X}_{N}:=\left\{  \psi_{j}\right\}  _{j=1}^{N}$
provides the "best possible approximation among the linear spaces of dimension
$N$" in the following sense: If we put
\begin{equation}
d_{N}\left(  K_{p}\right)  :=\inf_{X_{N}}\operatorname*{dist}\left(
X_{N},K_{p}\right)  \label{KolmogorovWidth}%
\end{equation}
the main result he proved in \cite{kolmogorov1936} says
\begin{equation}
d_{N}\left(  K_{p}\right)  =\operatorname*{dist}\left(  \widetilde{X}%
_{N},K_{p}\right)  . \label{KolmogorovWidth2}%
\end{equation}
Here we have used the notations, to be used also further,
\begin{align}
\operatorname*{dist}\left(  X,K_{p}\right)   &  :=\sup_{y\in K_{p}%
}\operatorname*{dist}\left(  X,y\right) \label{dist1}\\
\operatorname*{dist}\left(  X,y\right)   &  =\inf_{x\in X}\left\Vert
x-y\right\Vert . \label{dist2}%
\end{align}
Hence, by inequality (\ref{Jackson1dim}), equality (\ref{KolmogorovWidth2})
reads as
\begin{align*}
d_{N}\left(  K_{p}\right)   &  =\frac{1}{\sqrt{\lambda_{N+1}}}\qquad\text{for
}N\geq p\\
d_{N}\left(  K_{p}\right)   &  =\infty\qquad\qquad\text{for }N=0,1,...,p-1.
\end{align*}

\begin{definition}
The left quantity in (\ref{KolmogorovWidth}) is called \textbf{Kolmogorov
}$N-$\textbf{width}, while the best approximation space $\widetilde{X}_{N}$ is
called \textbf{extremal (optimal) subspace} (cf. this terminology in
\cite{tikhomirov}, \cite{lorentz}, \cite{pinkus}).
\end{definition}

Thus the \emph{main approach to the successful application of the theory of
widths} is based on a Jackson type theorem by which a special space
$\widetilde{X}_{N}$ is identified. Then one has to find, among which subspaces
$X_{N}$ is $\widetilde{X}_{N}$ the extremal subspace. Put in a different
perspective : one has to find as wide class of spaces $X_{N}$ as possible,
among which $\widetilde{X}_{N}$ is the extremal subspace.

Now let us consider the following set which is a \emph{natural multivariate
generalization} of the above set $K_{p}$ defined in (\ref{Kp}): For a bounded
domain $B$ in $\mathbb{R}^{n}$ we put (more generally than (\ref{Kpstar}))
\begin{equation}
K_{p}^{\ast}:=\left\{  u\in H^{2p}\left(  B\right)  :%
{\displaystyle\int_{B}}
\left\vert L_{2p}u\left(  x\right)  \right\vert ^{2}dx\leq1\right\}  ,
\label{KpstarGeneral}%
\end{equation}
where $L_{2p}$ is a \emph{strongly elliptic} operator in $B.$ Let us remark
that the \emph{Sobolev space} $H^{2p}\left(  B\right)  $ is the multivariate
version of the space of absolutely continuous functions on the interval with a
highest derivative in $L_{2}$ (as in (\ref{Kp})). An important feature of the
set $K_{p}^{\ast}$ is that it contains an infinite-dimensional subspace
\[
\left\{  u\in H^{2p}\left(  B\right)  :L_{2p}u\left(  x\right)  =0,\quad
\text{for }x\in B\right\}  .
\]
Hence, all Kolmogorov widths are equal to infinity, i.e.
\begin{equation}
d_{N}\left(  K_{p}^{\ast}\right)  =\infty\qquad\text{for }N\geq0
\label{dNisInfinity}%
\end{equation}
and \emph{no way} is seen to improve this if one remains within the
finite-dimensio\-nal setting.

The main purpose of the present paper is to find a proper setting in the
framework of the Polyharmonic Paradigm which generalizes the above results of Kolmogorov.

\section{A reminder on Elliptic Boundary Value \allowbreak Problems
\label{SreminderBVP}}

Let us specify the properties of the domains and the elliptic operators which
we will consider. In what follows we assume that the domain $D,$ the
differential operators and the boundary operators satisfy conditions for
\textbf{regular Elliptic BVP}. Namely, we give the following:

\begin{definition}
\label{Delliptic}We will say that the system of operators
\[
\left\{  A;B_{j},\ j=1,2,...,m\right\}
\]
forms a \textbf{regular Elliptic BVP in the domain }$D\subset\mathbb{R}^{n}$
if the following conditions hold:

1. The operator
\[
A\left(  x,D_{x}\right)  =%
{\displaystyle\sum_{\left\vert \alpha\right\vert ,\left\vert \beta\right\vert
\leq m}}
\left(  -1\right)  ^{\left\vert \alpha\right\vert }D^{\alpha}a_{\alpha\beta
}\left(  x\right)  D^{\beta}%
\]
is a differential operator with a principal part defined as
\[
A_{0}\left(  x,D_{x}\right)  =%
{\displaystyle\sum_{\left\vert \alpha\right\vert +\left\vert \beta\right\vert
=2m}}
\left(  -1\right)  ^{\left\vert \alpha\right\vert }a_{\alpha\beta}\left(
x\right)  D^{\alpha+\beta}.
\]
It is \textbf{uniformly} \textbf{strongly elliptic}, i.e. for every $x\in D$
holds
\[
c_{0}\left\vert \xi\right\vert ^{2m}\leq\left\vert A_{0}\left(  x,\xi\right)
\right\vert \leq c_{1}\left\vert \xi\right\vert ^{2m}\qquad\text{for all real
}\xi\in\mathbb{R}^{n}\setminus\left\{  0\right\}  .
\]

2. The domain $D$ is bounded and has a boundary $\partial D$ of the class
$C^{2m}.$

3. For every pair of linearly independent real vectors $\xi,$ $\eta$ and
$x\in\overline{D}$ the polynomial in $z,$ $A_{0}\left(  x,\xi+z\eta\right)  $
has exactly $m$ roots with positive imaginary parts.

4. The coefficients of $A$ are in $C^{\infty}\left(  \overline{D}\right)  .$
The boundary operators \allowbreak\ $B_{j}\left(  x,D\right)  =%
{\displaystyle\sum_{\left\vert \alpha\right\vert \leq m_{j}}}
b_{j,\alpha}\left(  x\right)  D^{\alpha}$ form a \textbf{normal system}, i.e.
their principal symbols are \textbf{non-characteristic}, i.e. satisfy
$B_{j,0}\left(  x,\xi\right)  =%
{\displaystyle\sum_{\left\vert \alpha\right\vert =m_{j}}}
b_{j,\alpha}\left(  x\right)  \xi^{\alpha}\neq0$ for every $x\in\partial D$
and $\xi\neq0,$ $\xi$ is normal to $\partial D$ at $x;$ they have pairwise
different orders $m_{j}$ which satisfy $m_{j}<2m$ for $1\leq j\leq m,$ and
their coefficients $b_{j,\alpha}$ belong to $C^{\infty}$ in $\partial D.$

5. At any point $x\in\partial D$ let $\nu$ denote the outward normal to
$\partial D$ at $x$ and let $\xi\neq0$ be a real vector in the tangent
hyperplane to $\partial D$ at $x.$ The polynomials in $z$ given by
$B_{j,0}\left(  x,\xi+z\nu\right)  $ are linearly independent modulo the
polynomial $\prod_{k=1}^{m}\left(  z-z_{k}^{+}\left(  \xi\right)  \right)  $
where $z_{k}^{+}\left(  \xi\right)  $ denote the roots of $A_{0}\left(
x,\xi+z\eta\right)  $ with positive imaginary parts.
\end{definition}

\begin{remark}
With minor differences the above definition is available in
\cite{lions-magenes} (conditions (i)-(iii) in chapter $2,$ section $5.1$); in
\cite{taylor} (sections $5.11$ and $5.12$); in \cite{hoermanderVol3} (chapter
$20$); in \cite{okbook} (section $23.2,$ p. $473$).
\end{remark}

Let us define a special system of boundary operators called \emph{Dirichlet}.
We put
\begin{align*}
B_{j}  &  =\left(  \frac{\partial}{\partial n}\right)  ^{j-1}\qquad\text{for
}j=1,2,...,p-1\\
S_{j}  &  =\left(  \frac{\partial}{\partial n}\right)  ^{p+j-1}\qquad\text{for
}j=1,2,...,p-1.
\end{align*}
Obviously,
\[
\operatorname*{ord}\left(  B_{j}\right)  =j-1,\qquad\operatorname*{ord}\left(
S_{j}\right)  =p+j-1.
\]
Let us denote by $L_{2p}^{\ast}$ the operator formally adjoint to the elliptic
operator $L_{2p}$. There exist boundary operators $C_{j},$ $T_{j},$ for
$j=1,2,...,p-1,$ such that
\[
\operatorname*{ord}\left(  T_{j}\right)  =2p-j,\qquad\operatorname*{ord}%
\left(  C_{j}\right)  =p-j
\]
and the following Green's formula holds:\
\begin{equation}%
{\displaystyle\int_{B}}
\left(  L_{2p}u\cdot v-u\cdot L_{2p}^{\ast}v\right)  dx=%
{\displaystyle\sum_{j=0}^{p-1}}
{\displaystyle\int_{\partial B}}
\left(  S_{j}u\cdot C_{j}v-B_{j}u\cdot T_{j}v\right)  d\sigma_{y};
\label{GreenGeneral}%
\end{equation}
here $\partial_{n}$ denotes the normal derivative to $\partial B,$ for
functions $u$ and $v$ in the classes of Sobolev, $u,v$ $\in$ $H^{2p}\left(
B\right)  $ (cf. \cite{lions-magenes}, Theorem $2.1$ in section $2.2,$ chapter
$2$, and Remark $2.2$ in section $2.3$).

For us the following eigenvalue problem will be important to consider for
$U\in H^{2p}\left(  B\right)  ,$ which is analogous to problem (\ref{eigen1}%
)-(\ref{eigen2}):%
\begin{align}
L_{2p}^{\ast}L_{2p}U\left(  x\right)   &  =\lambda U\left(  x\right)
\qquad\qquad\qquad\text{for }x\in B\label{eigen1Multi}\\
B_{j}L_{2p}U\left(  y\right)   &  =S_{j}L_{2p}U\left(  y\right)
=0,\qquad\text{for }y\in\partial B,\quad j=0,1,...,p-1 \label{eigen2Multi}%
\end{align}
where $\partial_{n}$ denotes the normal derivative at $y\in\partial B.$ It is
obvious that the operator $L_{2p}^{\ast}L_{2p}$ is formally self-adjoint,
however the BVP (\ref{eigen1Multi})-(\ref{eigen2Multi}) is not a nice one.
Since a direct reference seems not to be available, we provide its
consideration in the following theorem which is an analog to Proposition
\ref{PBVP}.

\begin{theorem}
\label{TExpansionBerezanskii} Let the operator $L_{2p}$ be uniformly strongly
elliptic in the domain $B.$ Then problem (\ref{eigen1Multi}%
)-(\ref{eigen2Multi}) has only real non-negative eigenvalues.

1. The eigenvalue $\lambda=0$ has infinite multiplicity with corresponding
eigenfunctions $\left\{  \psi_{j}^{\prime}\right\}  _{j=1}^{\infty}$ which
represent an orthonormal basis of the space of all solutions to the equation
$L_{2p}U\left(  x\right)  =0,$ for $x\in B.$

2. The positive eigenvalues are countably many and each has \textbf{finite
multiplicity}, and if we denote them by $\lambda_{j}$ ordered increasingly,
they satisfy $\lambda_{j}\longrightarrow\infty$ for $j\longrightarrow\infty.$

3. The orthonormalized eigenfunctions, corresponding to eigenvalues
$\lambda_{j}>0,$ will be denoted by $\left\{  \psi_{j}\right\}  _{j=1}%
^{\infty}.$ The set of functions $\left\{  \psi_{j}\right\}  _{j=1}^{\infty
}\bigcup\left\{  \psi_{j}^{\prime}\right\}  _{j=1}^{\infty}$ form a complete
orthonormal system in $L_{2}\left(  B\right)  .$
\end{theorem}

\begin{remark}
Problem (\ref{eigen1Multi})-(\ref{eigen2Multi}) is well known to be a
non-regular elliptic BVP, as well as non-coercive variational, cf.
\cite{agmon} ( p. $150$ ) and \cite{lions-magenes} (Remark $9.8$ in chapter
$2,$ section $9.6$, and section $9.8$ ).
\end{remark}

The proof is provided in the Appendix below, section \ref{Sappendix}.

\section{The principal axes of the ellipsoid $K_{p}^{\ast}$ and a Jackson type
theorem \label{SprincipalAxes}}

Here we will find the principal exes of the ellipsoid $K_{p}^{\ast}$ defined
as
\begin{equation}
K_{p}^{\ast}:=\left\{  u\in H^{2p}\left(  B\right)  :%
{\displaystyle\int_{B}}
\left\vert L_{2p}u\left(  x\right)  \right\vert ^{2}dx\leq1\right\}  ,
\label{KpstarGENERAL}%
\end{equation}
where $L_{2p}$ is a \emph{uniformly strongly elliptic} operator in $B.$

We prove the following theorem which generalizes Kolmogorov's one-dimensional
result from Proposition \ref{PKolmogorovJackson}, about the representation of
the ellipsoid $K_{p}$ in principal axes.

\begin{theorem}
\label{TprincipalAxesMultivariate}Let $f\in K_{p}^{\ast}.$ Then $f$ is
represented in a $L_{2}-$series as
\[
f\left(  x\right)  =%
{\displaystyle\sum_{j=1}^{\infty}}
f_{j}^{\prime}\psi_{j}^{\prime}\left(  x\right)  +%
{\displaystyle\sum_{j=1}^{\infty}}
f_{j}\psi_{j}\left(  x\right)  ,
\]
where by Theorem \ref{TExpansionBerezanskii} the eigenfunctions $\psi
_{j}^{\prime}$ satisfy $\Delta^{p}\psi_{j}^{\prime}\left(  x\right)  =0$ while
the eigenfunctions $\psi_{j}$ correspond to the eigenvalues $\lambda_{j}>0,$
and also
\begin{equation}%
{\displaystyle\sum_{j=1}^{\infty}}
\lambda_{j}f_{j}^{2}\leq1. \label{EllipsoidCondition}%
\end{equation}
Vice versa, every sequence $\left\{  f_{j}^{\prime}\right\}  _{j=1}^{\infty
}\bigcup\left\{  f_{j}\right\}  _{j=1}^{\infty}$ with $%
{\displaystyle\sum_{j=1}^{\infty}}
\left\vert f_{j}^{\prime}\right\vert ^{2}+%
{\displaystyle\sum_{j=1}^{\infty}}
\left\vert f_{j}\right\vert ^{2}<\infty$ and $%
{\displaystyle\sum_{j=1}^{\infty}}
\lambda_{j}f_{j}^{2}\leq1$ defines a function $f\in L_{2}\left(  B\right)  $
which is in $K_{p}^{\ast}.$
\end{theorem}

%

\proof
\textbf{(1)} According to Theorem \ref{TExpansionBerezanskii}, we know that
arbitrary $f\in L_{2}\left(  B\right)  $ is represented as
\begin{align*}
f\left(  x\right)   &  =%
{\displaystyle\sum_{j=1}^{\infty}}
f_{j}^{\prime}\psi_{j}^{\prime}\left(  x\right)  +%
{\displaystyle\sum_{j=1}^{\infty}}
f_{j}\psi_{j}\left(  x\right) \\
\left\Vert f\right\Vert _{L_{2}}^{2}  &  =%
{\displaystyle\sum_{j=1}^{\infty}}
\left\vert f_{j}^{\prime}\right\vert ^{2}+%
{\displaystyle\sum_{j=1}^{\infty}}
\left\vert f_{j}\right\vert ^{2}<\infty
\end{align*}
with convergence in the space $L_{2}\left(  B\right)  .$

\textbf{(2)} From the proof of Theorem \ref{TExpansionBerezanskii}, we know
that if we put
\[
\phi_{j}\left(  x\right)  =L_{2p}\psi_{j}\left(  x\right)  \qquad\text{for
}j\geq1,
\]
then the system of functions
\[
\frac{\phi_{j}\left(  x\right)  }{\sqrt{\lambda_{j}}}\qquad\text{for }j\geq1
\]
is orthonormal sequence which is complete in $L_{2}\left(  B\right)  .$

\textbf{(3)} We will prove now that if $f\in L_{2}\left(  B\right)  $ then
$f\in K_{p}^{\ast}$ iff
\[%
{\displaystyle\sum_{j=1}^{\infty}}
f_{j}^{2}\lambda_{j}\leq1.
\]
Indeed, for every $f\in H^{2p}\left(  B\right)  $ we have the expansion
$f\left(  x\right)  =%
{\displaystyle\sum_{j=1}^{\infty}}
f_{j}^{\prime}\psi_{j}^{\prime}\left(  x\right)  +%
{\displaystyle\sum_{j=1}^{\infty}}
f_{j}\psi_{j}\left(  x\right)  .$ We want to see that it is possible to
differentiate termwise this expansion, i.e.
\[
L_{2p}f\left(  x\right)  =%
{\displaystyle\sum_{j=1}^{\infty}}
f_{j}L_{2p}\psi_{j}\left(  x\right)  =%
{\displaystyle\sum_{j=1}^{\infty}}
f_{j}\phi_{j}\left(  x\right)
\]
Since $\left\{  \frac{\phi_{j}}{\sqrt{\lambda_{j}}}\right\}  _{j\geq1}$ is a
complete orthonormal basis of $L_{2}\left(  B\right)  $ it is sufficient to
see that
\[%
{\displaystyle\int_{B}}
L_{2p}f\left(  x\right)  \phi_{j}dx=%
{\displaystyle\int_{B}}
\left(
{\displaystyle\sum_{j=1}^{\infty}}
f_{j}L_{2p}\psi_{j}\left(  x\right)  \right)  \phi_{j}dx.
\]
Due to the boundary properties of $\phi_{j}$ and since $\phi_{j}=L_{2p}%
\psi_{j},$ we obtain
\[%
{\displaystyle\int_{B}}
L_{2p}f\left(  x\right)  \phi_{j}dx=%
{\displaystyle\int_{B}}
f\left(  x\right)  L_{2p}^{\ast}\phi_{j}dx=\lambda_{j}%
{\displaystyle\int_{B}}
f\psi_{j}dx=\lambda_{j}f_{j}.
\]
On the other hand
\[%
{\displaystyle\int_{B}}
\left(
{\displaystyle\sum_{k=1}^{\infty}}
f_{k}\phi_{k}\left(  x\right)  \right)  \phi_{j}dx=\lambda_{j}f_{j}.
\]
Hence
\[
L_{2p}f\left(  x\right)  =%
{\displaystyle\sum_{j=1}^{\infty}}
f_{j}L_{2p}\psi_{j}\left(  x\right)  =%
{\displaystyle\sum_{j=1}^{\infty}}
f_{j}\phi_{j}\left(  x\right)  =%
{\displaystyle\sum_{j=1}^{\infty}}
\sqrt{\lambda_{j}}f_{j}\frac{\phi_{j}\left(  x\right)  }{\sqrt{\lambda_{j}}}%
\]
and since $\left\{  \frac{\phi_{j}}{\sqrt{\lambda_{j}}}\right\}  _{j\geq1}$ is
an orthonormal system, it follows
\[
\left\Vert L_{2p}f\right\Vert _{L_{2}}^{2}=%
{\displaystyle\sum_{j=1}^{\infty}}
\lambda_{j}f_{j}^{2}.
\]
Thus if $f\in K_{p}$ it follows that $%
{\displaystyle\sum_{j=1}^{\infty}}
\lambda_{j}f_{j}^{2}\leq1.$

Now, assume vice versa, that $%
{\displaystyle\sum_{j=1}^{\infty}}
f_{j}^{2}\lambda_{j}\leq1$ holds together with $%
{\displaystyle\sum_{j=1}^{\infty}}
\left\vert f_{j}^{\prime}\right\vert ^{2}+%
{\displaystyle\sum_{j=1}^{\infty}}
\left\vert f_{j}\right\vert ^{2}<\infty$. We have to see that the function
\[
f\left(  x\right)  =%
{\displaystyle\sum_{j=1}^{\infty}}
f_{j}^{\prime}\psi_{j}^{\prime}\left(  x\right)  +%
{\displaystyle\sum_{j=1}^{\infty}}
f_{j}\psi_{j}\left(  x\right)
\]
belongs to the space $H^{2p}\left(  B\right)  .$ Based on the completeness and
orthonormality of the system $\left\{  \frac{\phi_{j}\left(  x\right)  }%
{\sqrt{\lambda_{j}}}\right\}  _{j=1}^{\infty}$ we may define the function
$g\in L_{2}$ by putting
\[
g\left(  x\right)  =%
{\displaystyle\sum_{j=1}^{\infty}}
\sqrt{\lambda_{j}}f_{j}\frac{\phi_{j}\left(  x\right)  }{\sqrt{\lambda_{j}}}=%
{\displaystyle\sum_{j=1}^{\infty}}
f_{j}\phi_{j}\left(  x\right)  ;
\]
it obviously satisfies $\left\Vert g\right\Vert _{L_{2}}\leq1.$

From the local solvability of elliptic equations (\cite{lions-magenes}) there
exists a function $F\in H^{2p}\left(  B\right)  $ which is a solution to
equation $L_{2p}F=g.$ Let its representation be
\[
F\left(  x\right)  =%
{\displaystyle\sum_{j=1}^{\infty}}
f_{j}^{\prime}\psi_{j}^{\prime}\left(  x\right)  +%
{\displaystyle\sum_{j=1}^{\infty}}
F_{j}\psi_{j}\left(  x\right)
\]
with some coefficients $F_{j}$ satisfying $%
{\displaystyle\sum_{j}}
\left\vert F_{j}\right\vert ^{2}<\infty.$ As above we obtain
\begin{align*}
\lambda_{j}%
{\displaystyle\int_{B}}
F\psi_{j}dx  &  =%
{\displaystyle\int_{B}}
FL_{2p}^{\ast}L_{2p}\psi_{j}dx=%
{\displaystyle\int_{B}}
L_{2p}F\cdot L_{2p}\psi_{j}dx\\
&  =%
{\displaystyle\int_{B}}
g\cdot\phi_{j}dx
\end{align*}
which implies $F_{j}=f_{j}.$ Hence, $F=f$ and $f\in H^{2p}\left(  B\right)  .$
This ends the proof.%

\endproof

We are able to prove finally a \emph{Jackson type} result as in Proposition
\ref{PKolmogorovJackson}.

\begin{theorem}
Let $N\geq1.$ Then for every $N\geq1$ and every $f\in K_{p}^{\ast}$ holds the
following estimate:
\[
\left\Vert f-%
{\displaystyle\sum_{j=1}^{\infty}}
f_{j}^{\prime}\psi_{j}^{\prime}\left(  x\right)  -%
{\displaystyle\sum_{j=1}^{N}}
f_{j}\psi_{j}\left(  x\right)  \right\Vert _{L_{2}}\leq\frac{1}{\sqrt
{\lambda_{N+1}}}.
\]

\end{theorem}

%

\proof
The proof follows directly. Indeed, due to the monotonicity of $\lambda_{j},$
and inequality (\ref{EllipsoidCondition}), we obtain
\[
\left\Vert f-%
{\displaystyle\sum_{j=1}^{\infty}}
f_{j}^{\prime}\psi_{j}^{\prime}\left(  x\right)  -%
{\displaystyle\sum_{j=1}^{N}}
f_{j}\psi_{j}\left(  x\right)  \right\Vert _{L_{2}}^{2}=%
{\displaystyle\sum_{j=N+1}^{\infty}}
f_{j}^{2}\leq\frac{1}{\lambda_{N+1}}%
{\displaystyle\sum_{j=N+1}^{\infty}}
f_{j}^{2}\lambda_{j}\leq\frac{1}{\lambda_{N+1}}.
\]
This ends the proof.%

\endproof

\section{Introducing the Hierarchy and Harmonic Widths \label{Shierarchy}}

In the present section we introduce the simplest representatives of the class
of domains having Harmonic Dimension $N,$ which are called \textbf{First Kind}
domains. They are piece-wise solutions to regular elliptic equations$.$

\begin{definition}
\label{Dhdimension} Let $D\subset\mathbb{R}^{n}$ be a bounded domain. For an
integer $M\geq1$ we say that the linear subspace $X_{M}\subset L_{2}\left(
D\right)  $ is of \textbf{First Kind} and has \textbf{Harmonic Dimension} $M,$
and write
\begin{equation}
\operatorname{hdim}\left(  X_{M}\right)  =M, \label{hdim}%
\end{equation}
if the following conditions are fulfilled:

1. There exists a finite number of domains $D_{j}$ with \textbf{piece-wise
smooth} boundaries $\partial D_{j}$ (which guarantees the validity of Green's
formula (\ref{GreenGeneral})), \allowbreak\ which are pairwise disjoint, i.e.
$D_{i}\bigcap D_{j}=\varnothing$ for $i\neq j,$ and such that we have the
\textbf{domain partition}
\begin{equation}
D=\bigcup_{j}D_{j}. \label{domainPartition}%
\end{equation}

2. We assume that for $j=1,2,...,M$ the \textbf{factorization operators}
$Q_{j}=Q_{j}\left(  x,D_{x}\right)  $ are \textbf{second order uniformly
strongly elliptic} which satisfy the maximum principle in the domain $D,$ and
the functions $\rho_{j}$ defined in $D$ are infinitely smooth and satisfy
\[
Z:=\bigcup_{k=1}^{M}\left\{  x\in D:\rho_{k}\left(  x\right)  =0\right\}
\subset\bigcup_{j}\partial D_{j}\setminus\partial D.
\]
Define the operator
\begin{equation}
P_{2M}\left(  x,D_{x}\right)  u\left(  x\right)  =Q_{M}\frac{1}{\rho_{M}%
}Q_{M-1}\frac{1}{\rho_{M-1}}\cdot\cdot\cdot Q_{1}\frac{1}{\rho_{1}}u\left(
x\right)  \label{P2Nproduct}%
\end{equation}
for the points $x\in D$ where it is correctly defined (out of the set $Z$ ).

We specify the interface conditions:\ Let us denote by $u_{i}=u_{|D_{i}}$ the
restriction of $u$ to $D_{i}.$ If for some indexes $i\neq j,\ $the
intersection $H:=\partial D_{i}\bigcap\partial D_{j}$ has nonempty interior in
the relative topology of $\partial D_{i}$ (hence also in $\partial D_{j}$)
then the following \textbf{interface conditions} hold on $H$ in the sense of
traces:
\begin{equation}
\left(  \frac{\partial}{\partial n_{x}}\right)  ^{k}u_{i}\left(  x\right)
=\left(  \frac{\partial}{\partial n_{x}}\right)  ^{k}u_{j}\left(  x\right)
\qquad\text{for }k=0,1,...,2M-1; \label{Interfaces}%
\end{equation}
here the vector $n_{x}$ denotes one of the normals at $x$ to the surface
$\partial D_{i}\bigcap\partial D_{j}.$

We define the space $X_{M}$ by putting
\begin{equation}
X_{M}=\left\{
\begin{array}
[c]{c}%
u\in H^{2M}\left(  D\right)  :P_{2M}u\left(  x\right)  =0,\quad\text{for }%
x\in\bigcup_{j=1}^{k}D_{j},\ \\
\text{and }u\text{ satisfies the interface conditions (\ref{Interfaces})}%
\end{array}
\right\}  . \label{XNDj}%
\end{equation}

\end{definition}

\begin{definition}
In the case of trivial partition (\ref{domainPartition}), i.e. $D=D_{1}$ where
we have functions $\rho_{j}$ free of zeros in $D$ we call the space $X_{M}$
\textbf{Multidimensional Chebyshev space}.
\end{definition}

Hence, by the above definitions, if $X_{M}\subset L_{2}\left(  D\right)  $ is
of First Kind and has Harmonic Dimension $M,$ then its restrictions to every
subdomain $D_{j}$ is a Multidimensional Chebyshev space, i.e. $X_{M}$ is a
piecewise Multidimensional Chebyshev space.

\begin{remark}
1. In \cite{kounchevSozopol} we considered the case of spaces $X_{M}$ of
Harmonic Dimension $M$ defined by a single elliptic operator $P_{2M}$ (i.e.
$P_{2M}=Q_{1}$) and a trivial partition of $D,$ i.e. $D=D_{1}.$

2. Let us comment on the interface conditions (\ref{Interfaces}) in Definition
\ref{Dhdimension}. Let us assume that we have an elliptic operator $P_{2M}$
with smooth coefficients defined on $D$ and that a non-trivial partition
$\bigcup D_{j}$ is given. Due to the piece-wise smoothness of the boundaries
$\partial D_{j}$ we may apply the Green formula, and from the interface
conditions (\ref{Interfaces}) it follows that "analytic continuation" is
possible, hence every function in $X_{M}$ is a solution to $P_{2M}u=0$ in the
whole domain $D$ (see similar result in \cite{okbook}, Lemma $20.10,$ and the
proof of Theorem $20.11$).

3. One may choose a different set of interface conditions which are equivalent
to (\ref{Interfaces}), see \cite{okbook} (Remark $20.12$), and
\cite{lions-magenes} (Lemma $2.1$ in chapter $2$).

4. The spaces $X_{M}$ defined in Definition \ref{Dhdimension} mimic in a
natural way the one-dimensional case: the operator $P_{2M}$ (\ref{P2Nproduct})
is similar to the operator (\ref{LNoperatorForm}) in Proposition
\ref{PstructureMarkovSystem}.

5. The operator $\left(  \prod_{j}\rho_{j}\right)  \times P_{2M}$ does not
have a singularity in the principal symbol but possibly only in the lower
order coefficients.
\end{remark}

Here is a simple non-trivial example to Definition \ref{Dhdimension}:
\begin{align*}
D_{1}  &  =\left\{  x:\left\vert x\right\vert <1\right\}  ,\quad
D_{2}=\left\{  x:1<\left\vert x\right\vert <2\right\} \\
D  &  =\left\{  x:\left\vert x\right\vert <2\right\} \\
P_{4}^{1}\left(  x;D_{x}\right)  u\left(  x\right)   &  =\Delta\frac
{1}{1-\left\vert x\right\vert }\Delta u\left(  x\right)  \qquad\text{for }x\in
D_{1}\\
P_{4}^{2}\left(  x;D_{x}\right)  u\left(  x\right)   &  =-\Delta\frac
{1}{1-\left\vert x\right\vert }\Delta u\left(  x\right)  \qquad\text{for }x\in
D_{2},
\end{align*}
where $\Delta$ is the Laplace operator. Typical elements of $X_{2}$ are the
functions $u$ which are obtained as solutions to
\[
\Delta u=\left(  1-\left\vert x\right\vert \right)  w\qquad\text{in }D,
\]
where $\Delta w=0$ in $D.$

The following result shows that we may construct a lot of solutions belonging
to the set $X_{M}$ of Definition \ref{Dhdimension}. We call these "direct solutions".

\begin{proposition}
\label{Pdirectsolutions} Let us define the boundary conditions $B_{k},$
$k=1,2,...,M,$ on $\partial D,$ by putting
\begin{align*}
B_{1}u  &  =\frac{1}{\rho_{1}}u\\
B_{k}u  &  =\frac{1}{\rho_{k}}Q_{k-1}B_{k-1}\qquad\text{for }k\geq2.
\end{align*}
Then the BVP
\begin{align}
P_{2M}u\left(  x\right)   &  =0\qquad\qquad\text{for }x\in D\label{BVP1}\\
B_{k}u\left(  y\right)   &  =h_{k}\left(  y\right)  \qquad\text{for }%
y\in\partial D,\ \text{and }k=1,2,...,M \label{BVP2}%
\end{align}
is solvable for arbitrary data $\left\{  h_{k}\right\}  _{k=1}^{M}$ from the
corresponding Sobolev spaces, i.e. $h_{k}\in H^{2M-ord\left(  B_{k}\right)
-1/2}\left(  \partial D\right)  ,$ and the solution has the maximal
regularity, i.e. $u\in H^{2M}\left(  D\right)  .$
\end{proposition}

%

\proof
For every $j$ with $1\leq j\leq M$ we consider the elliptic BVP of Dirichlet
\begin{align}
Q_{j}w  &  =f\qquad\quad\text{on }D\label{Qjw}\\
w  &  =h\qquad\text{on }\partial D. \label{Qjw2}%
\end{align}
By the assumptions on the operators $Q_{j}$ it is unique (by the maximum
principle) and we will denote it by $I_{j}\left(  f,h\right)  .$ It is easy to
see that the solution to $P_{2M}u\left(  x\right)  =0$ is obtained inductively
as
\begin{equation}
u=\rho_{1}I_{1}\left(  \cdot\cdot\cdot\rho_{M-1}I_{M-1}\left(  \rho_{M}%
I_{M}\left(  0;h_{M}\right)  ;h_{M-1}\right)  \cdot\cdot\cdot\right)  ,
\label{directSolution}%
\end{equation}
where $h_{j}$ are arbitrary boundary data.

For simplicity of notation let us assume that $k=2,$ i.e. $P=Q_{1}\frac
{1}{\rho_{1}}Q_{2}\frac{1}{\rho_{2}}.$ Then the boundary conditions satisfied
by $u$ are obtained from
\begin{align*}
Q_{2}w  &  =0\qquad\text{on }D\\
B_{2}w  &  =h_{2}\qquad\text{on }\partial D,
\end{align*}
and
\begin{align*}
Q_{1}\left(  \frac{1}{\rho_{1}}u\right)   &  =\rho_{2}w\qquad\text{on }D\\
B_{1}\frac{1}{\rho_{1}}u  &  =h_{1}\qquad\text{on }\partial D.
\end{align*}
Hence, we obtain
\[
B_{2}\left(  \frac{1}{\rho_{2}}Q_{1}\left(  \frac{1}{\rho_{1}}u\right)
\right)  =h_{2}.
\]
Thus we see that the system of boundary operators on $\partial D$ is
\[
B_{1}u=\frac{1}{\rho_{1}}u,\quad B_{2}u=\frac{1}{\rho_{2}}Q_{1}\left(
\frac{1}{\rho_{1}}u\right)
\]
and satisfies the conditions of Definition \ref{Delliptic}, for normal system
of boundary operators. We may proceed inductively to prove the statement for
arbitrary $k\geq3.$%

\endproof

\begin{remark}
1. Formula (\ref{directSolution}) for representing the solution of
$P_{2M}u\left(  x\right)  =0$ in Proposition \ref{Pdirectsolutions} coincides
with the solution in formulas (\ref{explicitECT1})-(\ref{explicitECT4}) in the
one-dimensional case.

2. One may prove that the set of "direct solutions" obtained in Proposition
\ref{Pdirectsolutions} is dense in the whole space $X_{M}$ defined in
Definition \ref{Dhdimension}, but we will not need this fact.
\end{remark}

The following fundamental theorem shows that, as in the one-dimensional case,
on arbitrary small sub-domain $G$ in $D$ with $G\bigcap\left(  \bigcup\partial
D_{j}\right)  =\varnothing,$ the space $X_{M}$ with $\operatorname{hdim}%
\left(  X_{M}\right)  =M$ has the same Harmonic Dimension $M.$ From a
different point of view, it shows that a theorem of Runge-Lax-Malgrange type
is true also for elliptic operators with singular coefficients of the type of
operators $P_{2M}$ considered in Definition \ref{Dhdimension}.

\begin{theorem}
\label{Tdensity} Let the First Kind space $X_{M}$ satisfy Definition
\ref{Dhdimension} with
\[
\operatorname{hdim}\left(  X_{M}\right)  =M.
\]
Assume that the elliptic operator $P_{2M}$ which corresponds to the space
$X_{M}$ has factorization operators $Q_{j}$ (from (\ref{P2Nproduct}))
satisfying condition $\left(  U\right)  _{s}$ for uniqueness in the Cauchy
problem in the small.\footnote{The differential operator $P$ satisfies
condition $\left(  U\right)  _{s}$ for uniqueness in the Cauchy problem in the
small in $G$ provided that if $G_{1}$ is a connected open subset of $G$ and
$u\in C^{r}\left(  G_{1}\right)  $ is a solution to $P^{\ast}u=0$ and $u$ is
zero on a non-emplty subset of $G_{1}$ then $u$ is identically zero. Elliptic
operators with analytic coefficients satisfy this property (cf.
\cite{bersJohnschechter}, part $II,$ chapter $1.4;$ \cite{browder}, p.
$402$).} Let $G$ be a compact subdomain in some $D_{j},$ i.e. $G\bigcap\left(
\bigcup\partial D_{j}\right)  =\varnothing.$ Then the set of "direct
solutions" considered in Proposition \ref{Pdirectsolutions} is dense in
$L_{2}\left(  G\right)  $ in the space
\[
\left\{  u\in H^{2M}\left(  G\right)  :P_{2M}u=0\quad\text{in }G\right\}  .
\]

\end{theorem}

%

\proof
For simplicity of notations we assume that for the elliptic operator $P_{2M}$
associated with $X_{M},$ by Definition \ref{Dhdimension}, we have only two
factorizing operators $Q_{1}$ and $Q_{2},$ i.e. $P_{2M}u=Q_{2}\frac{1}%
{\rho_{2}}Q_{1}\left(  \frac{1}{\rho_{1}}u\right)  .$

Let us take a solution $u\in H^{2M}\left(  G\right)  $ to $P_{2M}u=0$ in $G.$
We use the solution of formula (\ref{directSolution})
\[
u=\rho_{1}I_{1}\left(  \rho_{2}I_{2}\left(  0;h_{2}\right)  ;h_{1}\right)  ,
\]
where the boundary data $h_{1}$ and $h_{2}$ are arbitrary in proper Sobolev
spaces. By the approximation theorem of Runge-Lax-Malgrange type (cf.
\cite{browder}, Theorem $4,$ and references there), which uses essentially
property $\left(  U\right)  _{s}$ of operator $Q_{2},$ we obtain a function
$w_{\varepsilon}$ which is a solution to $Q_{2}w_{\varepsilon}=0$ in $D$ and
such that
\[
\left\Vert I_{2}\left(  0;h_{2}\right)  -w_{\varepsilon}\right\Vert
_{L_{2}\left(  G\right)  }<\varepsilon.
\]
Next we apply the same approximation argument but with non-zero right-hand
side $\rho_{2}I_{2}\left(  0;h_{2}\right)  $ (cf. \cite{browder2}) to prove
the existence of a function $v_{\varepsilon}$ such that
\[
\left\Vert I_{1}\left(  \rho_{2}I_{2}\left(  0;h_{2}\right)  ;h_{1}\right)
-v_{\varepsilon}\right\Vert <C\varepsilon
\]
for some constant $C>0,$ where the constant $C$ depends on the functions
$\rho_{j}.$ Thus we obtain the function
\[
u_{\varepsilon}=\rho_{1}v_{\varepsilon}%
\]
which satisfies
\[
\left\Vert u_{\varepsilon}-u\right\Vert _{L_{2}\left(  G\right)  }%
<C_{1}\varepsilon,
\]
and is a "direct solution" in the sense of Proposition \ref{Pdirectsolutions}.%

\endproof

The following theorem studies the orthogonal complement $X_{N}\ominus X_{M}$
of two First Kind spaces where $M<N.$ While we will not need the whole
generality of the result proved, the proof shows that $X_{N}\ominus X_{M}$ has
at least $\operatorname{hdim}$ equal to $N-M.$

\begin{theorem}
\label{TtransversalSpaces}Let $M<N$ and the First Kind spaces $X_{M},$ $X_{N}$
satisfy Definition \ref{Dhdimension} with
\[
\operatorname{hdim}\left(  X_{M}\right)  =M,\quad\operatorname{hdim}\left(
X_{N}\right)  =N.
\]
Assume that the elliptic operator $P_{2N}^{\prime},$ which is associated with
the space $X_{N},$ has (by (\ref{P2Nproduct})) factorization operators $Q_{j}$
satisfying condition $\left(  U\right)  _{s}$ for uniqueness in the Cauchy
problem in the small (as in Theorem \ref{Tdensity}). Then the space
$Y=X_{N}\setminus X_{M}$ is infinite-dimensional.
\end{theorem}

%

\proof
\textbf{(1)} Let, by Definition \ref{Dhdimension}, the partition $\bigcup
D_{j}$ and the operator $P_{2M}$ correspond to $X_{M},$ while the partition
$\bigcup D_{j}^{\prime}$ and the operator $P_{2N}^{\prime}$ correspond to
$X_{N}.$ Assume that $D_{1}\bigcap D_{1}^{\prime}\neq\varnothing.$ Then we
will choose a subdomain $G$ which is compactly supported in $D_{1}\bigcap
D_{1}^{\prime}.$

Further we will fix our attention to the subdomain $G$ where both operators
$P_{2M}$ and $P_{2N}^{\prime}$ are uniformly strongly elliptic and will
construct a subset of $X_{N}\ominus X_{M}$ restricted to the domain $G.$ Let
us be more precise: If we denote by
\begin{equation}
X_{N}^{G}:=\left\{  u:H^{2N}\left(  G\right)  :P_{2N}^{\prime}u=0\quad\text{in
}G\right\}  \label{XNG}%
\end{equation}
then we will construct an infinite-dimensional subspace of $X_{N}%
^{G}\circleddash X_{M}^{G}.$

\textbf{(2)} For the \emph{uniformly strongly elliptic} operator $P_{2M}$ on
the domain $G$ we choose the Dirichlet system of boundary operators
$B_{j}=\frac{\partial^{j-1}}{\partial n^{j-1}},$ for $j\geq1,$ which are
iterates of the normal derivative $\frac{\partial}{\partial n}$ on the
boundary $\partial G.$ As already mentioned the system of operators $\left\{
P_{2M};\frac{\partial^{j}}{\partial n^{j}}:j=0,1,...,M-1\right\}  $ on $G$
forms a \emph{regular Elliptic BVP }(this is the Dirichlet Elliptic BVP for
the operator $P_{2M}$) (cf. \cite{lions-magenes}, Remark $1.3$ in section
$1.4$, chapter $2$).

We complete the system $\left\{  B_{j}\right\}  _{j=1}^{M}$ by the system of
boundary operators $S_{j}=\frac{\partial^{M-1+j}}{\partial n^{M-1+j}}$ for
$j=1,2,...M.$ Hence, the system composed \allowbreak\ $\left\{  B_{j}\right\}
_{j=1}^{M}\bigcup\left\{  S_{j}\right\}  _{j=1}^{M}$ is a \emph{Dirichlet
system} of order $2M$ (cf. \cite{lions-magenes}, Definition $2.1$ and Theorem
$2.1$ in section $2.2,$ chapter $2$). Further, by \cite{lions-magenes}
(Theorem $2.1$), there exists a unique Dirichlet system of order $2M$ of
boundary operators $\left\{  C_{j},T_{j}\right\}  _{j=1}^{M}$ which is
uniquely determined as the adjoint to the system $\left\{  B_{j}%
,S_{j}\right\}  _{j=1}^{M},$ and the Green formula (\ref{GreenGeneral}) holds
on the domain $G.$ We will use this below.

\textbf{(3)} In the domain $G$ we consider the elliptic operator
$P_{2N}^{\prime}P_{2M}^{\ast}.$ As a product of two \emph{ strongly elliptic}
operators it is such again. By a standard construction cited above (cf.
\cite{lions-magenes}, Theorem $2.1,$ section $2.2,$ chapter $2$), we may
complete the Dirichlet system of operators $\left\{  B_{j},S_{j}\right\}
_{j=1}^{M}$ with $N-M$ boundary operators $R_{j}=\frac{\partial^{2M-1+j}%
}{\partial n^{2M-1+j}},$ $j=1,2,...,N-M.$ Again by the above cited theorem,
the Dirichlet system of boundary operators
\[
\left\{  B_{j},S_{j}\right\}  _{j=1}^{M}%
{\textstyle\bigcup}
\left\{  R_{j}\right\}  _{j=1}^{N-M}%
\]
\emph{covers} the operator $P_{2N}^{\prime}P_{2M}^{\ast}.$ Finally, we
consider the solutions \allowbreak\ $g\in H^{2N+2M}\left(  G\right)  $ to the
following Elliptic BVP:
\begin{align}
P_{2N}^{\prime}P_{2M}^{\ast}g\left(  x\right)   &  =0\qquad\qquad
\qquad\ \ \text{for }x\in G\label{gorthogonal1}\\
B_{j}g\left(  y\right)   &  =S_{j}g\left(  y\right)  =0\qquad\text{for
}j=0,1,...,N-1,\text{ for }y\in\partial G\label{gorthogonal2}\\
R_{j}g\left(  y\right)   &  =h_{j}\left(  y\right)  \qquad\qquad\ \ \text{for
}j=1,2,...,N-M,\text{ for }y\in\partial G. \label{gorthogonal3}%
\end{align}
We may apply a classical result \cite{lions-magenes} (the existence Theorem
$5.2$ and Theorem $5.3$ in chapter $2$), to the solvability of problem
(\ref{gorthogonal1})-(\ref{gorthogonal3}) in the space $H^{2M+2N}\left(
G\right)  .$

\textbf{(4) }Let us check the properties of the function $P_{2M}^{\ast}g$
where $g$ satisfies (\ref{gorthogonal1})-(\ref{gorthogonal3}). First of all,
it is clear from (\ref{gorthogonal1}) that $P_{2M}^{\ast}g\in X_{N}^{G}$ where
we have used the notation (\ref{XNG}).

By Green's formula (\ref{GreenGeneral}), applied for the operator $P_{2M}$ and
for $u=g$ we obtain
\[%
{\displaystyle\int_{G}}
P_{2M}^{\ast}g\cdot vdx=0\qquad\text{for all }v\text{ with }P_{2M}v=0
\]
which implies that the function $P_{2M}^{\ast}g$ satisfies $P_{2M}^{\ast
}g\perp X_{M}^{G}$ ($X_{M}^{G}$ defined as (\ref{XNG})).

By the general existence theorem for Elliptic BVP used already above (cf.
\cite{lions-magenes}, Theorem $5.3$, the Fredholmness property), we know that
a solution $g$ to problem (\ref{gorthogonal1})-(\ref{gorthogonal3}) exists for
those boundary data $\left\{  h_{j}\right\}  _{j=1}^{N-M}$ which satisfy only
a finite number of linear conditions (cf. \cite{lions-magenes}, conditions
(5.18)); these are determined by the solutions to the homogeneous adjoint
Elliptic BVP. Hence, it follows that the space $Y_{N-M}^{G}$ of the functions
$P_{2M}^{\ast}g$ where $g$ is a solution to (\ref{gorthogonal1}%
)-(\ref{gorthogonal3}) is infinite-dimensional.

\textbf{(5) }Let us construct a subspace of $X_{N}\setminus X_{M}$ which is
infinite-dimensional. We use the obvious inclusion $X_{N|G}\subset X_{N}^{G},$
$X_{M|G}\subset X_{M}^{G},$ where for a space of functions $Y\subset
L_{2}\left(  B\right)  $ the space $\ Y_{|G}$ consists of the restrictions of
the elements of $Y$ to the domain $G.$

First of all, we find an orthonormal basis $\left\{  v_{j}\right\}  _{j\geq1}$
in the infinite-dimensio\-nal space $Y_{N-M}^{G}$ (where the norm is
$\left\Vert \cdot\right\Vert _{L_{2}\left(  G\right)  }$ ); by the
Gram-Schmidt orthonormalization we obtain functions $g_{j}$ such that
$v_{j}=P_{2M}^{\ast}g_{j}$ for $j\geq1.$

Let us put $\varepsilon_{j}=\frac{1}{2^{j-1}}$ and use the density Theorem
\ref{Tdensity} to choose $u_{j}\in H^{2N}\left(  D\right)  $ with
\[
\left\Vert u_{j}-v_{j}\right\Vert _{L_{2}\left(  G\right)  }\leq
\varepsilon_{j}\qquad\text{for }j\geq1.
\]
The orthogonality of $v_{j}$ to $X_{M}^{G}$ infers $\operatorname*{dist}%
\left(  u_{j|G},X_{M}^{G}\right)  \geq1-\varepsilon_{j}$ in the $L_{2}\left(
G\right)  $ norm. Hence, $\operatorname*{dist}\left(  u_{j},X_{M}\right)
\geq1-\varepsilon_{j}$ in the $L_{2}\left(  D\right)  $ norm, hence
$u_{j}\notin X_{M}.$

Let us see that for every choice of the constants $\alpha_{j}$ holds
\[%
{\displaystyle\sum_{j=1}^{N-1}}
\alpha_{j}u_{j}\neq u_{N}.
\]
Indeed, by the triangle inequality for the norm $\left\Vert \cdot\right\Vert
_{L_{2}\left(  G\right)  }$ it follows
\begin{align*}
1+%
{\displaystyle\sum_{j=1}^{N-1}}
\left\vert \alpha_{j}\right\vert ^{2}  &  =\left\Vert v_{N}-%
{\displaystyle\sum_{j=1}^{N-1}}
\alpha_{j}v_{j}\right\Vert \\
&  =\left\Vert v_{N}-u_{N}+u_{N}-%
{\displaystyle\sum_{j=1}^{N-1}}
\alpha_{j}u_{j}+%
{\displaystyle\sum_{j=1}^{N-1}}
\alpha_{j}u_{j}-%
{\displaystyle\sum_{j=1}^{N-1}}
\alpha_{j}v_{j}\right\Vert \\
&  \leq\varepsilon_{N}+\left\Vert u_{N}-%
{\displaystyle\sum_{j=1}^{N-1}}
\alpha_{j}u_{j}\right\Vert +%
{\displaystyle\sum_{j=1}^{N-1}}
\left\vert \alpha_{j}\right\vert \varepsilon_{j}%
\end{align*}
or
\[
1-\varepsilon_{N}+%
{\displaystyle\sum_{j=1}^{N-1}}
\left(  \left\vert \alpha_{j}\right\vert ^{2}-\left\vert \alpha_{j}\right\vert
\varepsilon_{j}\right)  \leq\left\Vert u_{N}-%
{\displaystyle\sum_{j=1}^{N-1}}
\alpha_{j}u_{j}\right\Vert .
\]
Obviously
\[
1-\varepsilon_{N}+%
{\displaystyle\sum_{j=1}^{N-1}}
\left(  \frac{\varepsilon_{j}^{2}}{4}-\frac{\varepsilon_{j}}{2}\varepsilon
_{j}\right)  \leq1-\varepsilon_{N}+%
{\displaystyle\sum_{j=1}^{N-1}}
\left(  \left\vert \alpha_{j}\right\vert ^{2}-\left\vert \alpha_{j}\right\vert
\varepsilon_{j}\right)
\]
and since the left-hand side always exceeds $1/4$, this ends the proof that
the system of functions $\left\{  u_{j|G}\right\}  _{j\geq1}$ is linearly
independent. Hence, the system $\left\{  u_{j}\right\}  _{j\geq1}$ is linearly
independent in the whole domain $D.$

As noted above $u_{j}\notin X_{M},$ hence $\operatorname*{span}\left\{
u_{j}\right\}  _{j\geq1}$ is the infinite-dimensional space we sought. The
proof is finished.%

\endproof

We have the following prototype of Theorem \ref{TtransversalSpaces}, proved in
\cite{kounchevSozopol}.

\begin{corollary}
\label{Ctransversal}Let $M<N$ and $X_{M},$ $X_{N}$ satisfy Definition
\ref{Dhdimension} with
\[
\operatorname{hdim}\left(  X_{M}\right)  =M,\quad\operatorname{hdim}\left(
X_{N}\right)  =N.
\]
Assume that the differential operators $P_{2M}$ and $P_{2N}^{\prime},$
associated with $X_{M}$ and $X_{N},$ have trivial factorization operators by
the definition (\ref{P2Nproduct}), and trivial domain partitions $D=D_{1}$ and
$D=D_{1}^{\prime}$ by (\ref{domainPartition}). Then the space of solutions of
the Elliptic BVP (\ref{gorthogonal1})-(\ref{gorthogonal3}) where $G=D$ is a
subspace of the space
\[
Y=X_{N}\ominus X_{M}.
\]

\end{corollary}

The proof may be derived from the proof of Theorem \ref{TtransversalSpaces}
where we have put $G=D.$ Note that we \textbf{do not need} the $\left(
U\right)  _{s}$ condition for the operator $P_{2N}^{\prime}.$ Hence, strictly
speaking, Corollary \ref{Ctransversal} is not a special case of Theorem
\ref{TtransversalSpaces}.

Now we provide a generalization of Kolmogorov's notion of width from formula
(\ref{KolmogorovWidth}); without restricting the generality we assume that we
work only with symmetric subsets.

\begin{definition}
\label{Dwidth} Let $A$ be a centrally symmetric subset in $L_{2}\left(
B\right)  .$ For fixed integers $M\geq1$ and $N\geq0$ we define the
corresponding \textbf{Harmonic Width} by putting
\[
\operatorname*{hd}\nolimits_{M,N}\left(  K\right)  :=\inf_{X_{M},F_{N}%
}\operatorname*{dist}\left(  X_{M}%
{\textstyle\bigoplus}
F_{N},A\right)  ,
\]
where $\inf_{X_{M},F_{N}}$ is taken over all spaces $X_{M},F_{N}\subset
C^{\infty}\left(  B\right)  $ with
\begin{align*}
\operatorname{hdim}\left(  X_{M}\right)   &  =M\\
\dim\left(  F_{N}\right)   &  =N.
\end{align*}

\end{definition}

\section{Generalization of Kolmogorov's result \allowbreak\ about widths
\label{Swidths}}

Next we prove results which are analogs to the original Kolmogorov's results
about widths in (\ref{KolmogorovWidth2}).

We denote by $F_{N}$ a \textbf{finite-dimensional} subspace of $L_{2}\left(
B\right)  $ of dimension $N.$ We denote the special subspaces for an elliptic
operator $P_{2p}=L_{2p}$ by
\begin{equation}
\widetilde{X}_{p}:=\left\{  u\in H^{2p}\left(  B\right)  :L_{2p}u\left(
x\right)  =0,\quad\text{for }x\in B\right\}  , \label{StildeM}%
\end{equation}
and the special finite-dimensional subspaces
\begin{equation}
\widetilde{F}_{N}:=\left\{  \psi_{j}:j\leq N\right\}  _{lin} \label{FtildeN}%
\end{equation}
where $\psi_{j}$ are the eigenfunctions from Theorem
\ref{TExpansionBerezanskii}.

\begin{theorem}
\label{TKolmogorovMultivariate} Let $K_{p}^{\ast}$ be the set defined in
(\ref{KpstarGENERAL}) as
\[
K_{p}^{\ast}:=\left\{  u\in H^{2p}\left(  B\right)  :%
{\displaystyle\int_{B}}
\left\vert L_{2p}u\left(  x\right)  \right\vert ^{2}dx\leq1\right\}  ,
\]
with a constant coefficient operator $L_{2p}$ which is uniformly strongly
elliptic in the domain $B.$ Let $X_{M}$ be a First Kind subspace of
$L_{2}\left(  B\right)  $ of Harmonic Dimension $M,$ according to Definition
\ref{Dhdimension}, i.e.
\[
\operatorname{hdim}\left(  X_{M}\right)  =M,
\]
and let $N\geq0$ be arbitrary.

1. If $M<p$ then
\[
\operatorname*{dist}\left(  X_{M}%
{\textstyle\bigoplus}
F_{N},K_{p}^{\ast}\right)  =\infty.
\]
Hence,
\[
\inf_{X_{M},F_{N}}\operatorname*{dist}\left(  X_{M}%
{\textstyle\bigoplus}
F_{N},K_{p}^{\ast}\right)  =\infty
\]
or equivalently,
\[
\operatorname*{hd}\nolimits_{M,N}\left(  K_{p}^{\ast}\right)  =\infty.
\]

2. If $M=p$ then
\[
\inf_{X_{p},F_{N}}\operatorname*{dist}\left(  X_{p}%
{\textstyle\bigoplus}
F_{N},K_{p}^{\ast}\right)  =\operatorname*{dist}\left(  \widetilde{X}_{p}%
{\textstyle\bigoplus}
\widetilde{F}_{N},K_{p}^{\ast}\right)  ,
\]
i.e.
\[
\operatorname*{hd}\nolimits_{p,N}\left(  K_{p}^{\ast}\right)
=\operatorname*{dist}\left(  \widetilde{X}_{p}%
{\textstyle\bigoplus}
\widetilde{F}_{N},K_{p}^{\ast}\right)  .
\]

\end{theorem}

\begin{remark}
In both cases we see that the special spaces $\widetilde{X}_{M}%
{\textstyle\bigoplus}
\widetilde{F}_{N}$ are extremizers among the large class of spaces $X_{M}%
{\textstyle\bigoplus}
F_{N}.$
\end{remark}

%

\proof
1. If we assume that $X_{M}$ and $\widetilde{X}_{p}$ are transversal the proof
is clear since $\widetilde{X}_{p}\subset K_{p}^{\ast}$ and there will be an
infinite-dimensional subspace in $\widetilde{X}_{p}\subset K_{p}^{\ast}$
containing at least one infinite axis with direction $f\in\widetilde{X}%
_{p}\setminus X_{M},$ such that
\[
\operatorname*{dist}\left(  X_{M}%
{\textstyle\bigoplus}
F_{N},f\right)  >0
\]
which implies
\[
\operatorname*{dist}\left(  X_{M}%
{\textstyle\bigoplus}
F_{N},K_{p}^{\ast}\right)  =\infty.
\]
If they are not transversal we remind that operators with analytic
coefficients satisfy the $\left(  U\right)  _{s}$ condition, and we may apply
Lemma \ref{LtransversalSpaces}.

2. For proving the second item, let us first note that $\widetilde{X}%
_{p}\subset X_{p}%
{\textstyle\bigoplus}
F_{N}.$ Indeed, since $\widetilde{X}_{p}\subset K_{p}^{\ast}$ the violation of
$\widetilde{X}_{p}\subset X_{p}%
{\textstyle\bigoplus}
F_{N}$ would imply that there exists an infinite axis $f$ in $K_{p}^{\ast}$
not contained in $X_{p}%
{\textstyle\bigoplus}
F_{N}$ which would immediately give
\[
\operatorname*{dist}\left(  X_{p}%
{\textstyle\bigoplus}
F_{N},K_{p}^{\ast}\right)  =\infty.
\]
Using the notations of Definition \ref{Dhdimension}, there exists a finite
cover $\bigcup D_{j}=B,$ and by Lemma \ref{LellipticCanonical} (applied for
$M=N=p$ ) it follows that on every subdomain $D_{j}$ holds $P_{2p}^{j}%
=C_{j}\left(  x\right)  L_{2p}$ for some function $C_{j}\left(  x\right)  .$
Thus we see that every $u\in X_{p}$ is a piecewise solution of $L_{2p}u=0$ on
$B,$ satisfying the interface conditions (\ref{Interfaces}) in Definition
\ref{Dhdimension}. Here we use an uniqueness theorem for "analytic
continuation" across the boundary argument (proved directly by Green's formula
(\ref{GreenGeneral}) as in \cite{okbook}, Lemma $20.10$ and the proof of
Theorem $20.11,$ p. $422$) that $u\in\widetilde{X}_{p},$ hence $X_{p}%
=\widetilde{X}_{p}.$

Further we follow the usual way as in \cite{lorentz} to see that
$\widetilde{F}_{N}$ is extremal among all finite-dimensional spaces $F_{N},$
i.e.
\[
\inf_{F_{N}}\operatorname*{dist}\left(  \widetilde{X}_{p}%
{\textstyle\bigoplus}
F_{N},K_{p}^{\ast}\right)  =\operatorname*{dist}\left(  \widetilde{X}_{p}%
{\textstyle\bigoplus}
\widetilde{F}_{N},K_{p}^{\ast}\right)  .
\]
This ends the proof.%

\endproof

We prove the following fundamental result which shows the mutual position of
two subspaces:

\begin{lemma}
\label{LtransversalSpaces}Assume the conditions of Theorem
\ref{TtransversalSpaces}. Let the integer $M_{1}\geq0.$ Then%
\[
\operatorname*{dist}\left(  X_{M}%
{\textstyle\bigoplus}
F_{M_{1}},X_{N}\right)  =\infty.
\]

\end{lemma}

The proof follows directly from Theorem \ref{TtransversalSpaces} since a
finite-dimensional subspace $F_{M_{1}}$ would not disturb the arguments there.

We obtain immediately the following result.

\begin{corollary}
Let us denote by $U_{N+1}$ the unit ball in $X_{N+1}$ in the $L_{2}\left(
B\right)  $ norm. Then
\[
\operatorname*{dist}\left(  X_{N},U_{N+1}\right)  =1.
\]

\end{corollary}

\begin{remark}
Lemma \ref{LtransversalSpaces} and especially the above Corollary may be
considered as a generalization in our setting of a theorem of Gohberg-Krein of
$1957$ (cf. \cite{lorentz}, Theorem $2$ on p. $137$ ) in a Hilbert space.
\end{remark}

We need the following intuitive result which is however not trivial.

\begin{lemma}
\label{LellipticCanonical}Let for the strongly elliptic differential operators
\allowbreak\ $L_{2N}=P_{2N}\left(  x;D_{x}\right)  $ and $P_{2M}=P_{2M}\left(
x;D_{x}\right)  $ of orders respectively $2N\leq2M$ in the domain $B,$ the
following inclusion hold
\[
X_{N}%
{\textstyle\bigcap}
H^{2M}\left(  B\right)  \subset X_{M}\setminus F,
\]
or
\begin{align*}
&  \left\{  u\in H^{2M}\left(  B\right)  :L_{2N}u\left(  x\right)  =0,\quad
x\in B\right\}  \subset\\
&  \subset\left\{  u\in H^{2M}\left(  B\right)  :P_{2M}u\left(  x\right)
=0,\quad x\in B\right\}  \setminus F,
\end{align*}
where $F\subset L_{2}\left(  B\right)  $ is a \textbf{finite-dimensional}
subspace of $L_{2}\left(  B\right)  .$ Then
\begin{equation}
P_{2M}\left(  x,D_{x}\right)  =P_{2M-2N}^{\prime}\left(  x,D_{x}\right)
L_{2N}\left(  x,D_{x}\right)  \label{P2N=cx}%
\end{equation}
for some strongly elliptic differential operator $P_{2M-2N}^{\prime}$ of order
$2M-2N.$
\end{lemma}

%

\proof
It is clear that the arguments for proving equality (\ref{P2N=cx}) are purely
local, and it suffices to consider only $x_{0}=0,$ or we assume that the
operator $L_{2N}$ has constant coefficients.

First, we assume that the polynomial $L_{2N}\left(  \zeta\right)  $ is
\emph{irreducible}. Then we consider the roots of the equation
\begin{equation}
L_{2N}\left(  \zeta\right)  =0\qquad\text{for }\zeta\in\mathbb{C}^{n}.
\label{L2Nzeta=0}%
\end{equation}
If $\zeta$ is a solution to (\ref{L2Nzeta=0}) then the function $v\left(
x\right)  =\exp\left(  \left\langle \zeta,x\right\rangle \right)  $ is a
solution to equation $L_{2N}v=0$ in the whole space. Hence
\[
P_{2M}v=P_{2M}\left(  x_{0};D_{x}\right)  v\left(  x_{0}\right)
=P_{2M}\left(  x_{0};\zeta\right)  v\left(  x_{0}\right)  =0,
\]
and by a well-known result on division of polynomials in algebra \cite{walker}
(Theorem $9.7,$ p. $26$), the statement of the theorem follows.

Now let us assume that $L_{2N}$ is reducible and decomposed in two irreducible
factors $L_{2N}=Q_{2}Q_{1},$ which may be equal. Obviously, both polynomials
$Q_{1}$ and $Q_{2}$ are uniformly strongly elliptic. Since the solutions to
$Q_{1}u=0$ are also solutions to $L_{2N}$ it follows by the above that
\[
P_{2M}\left(  x,D_{x}\right)  =P_{2M-2N_{1}}^{\prime}\left(  x,D_{x}\right)
Q_{1}\left(  D_{x}\right)
\]
where $2N_{1}$ is the order of the operator $Q_{1}.$ Further, following the
standard arguments in \cite{lions-magenes}, by the uniform strong ellipticity
of the operator $Q_{1},$ for every $\zeta\in\mathbb{C}^{n},$ and for arbitrary
$s\geq2N_{1},$ there exists a solution $u\in H^{s}\left(  B\right)  $ to
equation
\[
Q_{1}u_{\zeta}\left(  x\right)  =e^{\left\langle \zeta,x\right\rangle }%
\qquad\text{for }x\in B.
\]
Let $\zeta\in\mathbb{C}^{n}$ be a solution to equation $Q_{2}\left(
\zeta\right)  =0.$ Obviously,
\[
L_{2N}u_{\zeta}=0
\]
hence, by the above it follows
\[
P_{2M}\left(  x,D_{x}\right)  u_{\zeta}=P_{2M-2N_{1}}^{\prime}\left(
x,D_{x}\right)  Q_{1}\left(  D_{x}\right)  u_{\zeta}=P_{2M-2N_{1}}^{\prime
}\left(  x,\zeta\right)  =0.
\]
It follows that $P_{2M-2N_{1}}^{\prime}\left(  x_{0},\zeta\right)  =0.$ We
proceed inductively if $L_{2N}$ has more than two irreducible factors.%

\endproof

\section{Second Kind spaces of Harmonic Dimension $N$ and widths
\label{SsecondKind}}

In order to make things more transparent, in Definition \ref{Dhdimension} we
avoided the maximal generality of the notions and considered only First Kind
spaces of Harmonic Dimension $N.$ Let us explain by analogy with the
one-dimensional case how do the "Second Kind" spaces of Harmonic Dimension $N$ appear.

In the one-dimensional case, if we have a finite-dimensional subspace
$X_{N}\subset C^{N}\left(  I\right)  $ then for a point $x_{0}\in I$ the
space
\[
Y:=\left\{  u\in X_{N}:u\left(  x_{0}\right)  =0\right\}
\]
is an $\left(  N-1\right)  -$dimensional subspace. We would like that our
notion of Harmonic Dimension $N$ behave in a similar way. For example, if
$X_{N}$ is defined as a set of solutions of an elliptic operator $P_{2N}$ by
\[
X_{N}:=\left\{  u\in H^{2N}\left(  B\right)  :P_{2N}u=0\quad\text{in
}B\right\}
\]
then it is natural to expect that the space
\[
Y:=\left\{  u\in X_{N}:u=0\quad\text{on }\partial B\right\}
\]
has Harmonic Dimension $N-1.$ A simple example is the space
\[
Y=\left\{  u\in H^{4}\left(  B\right)  :\Delta^{2}u=0\quad\text{in
}B,\ u=0\quad\text{on }\partial B\right\}  .
\]

On the other hand, it is Theorem \ref{TtransversalSpaces} and Corollary
\ref{Ctransversal} above which show that such Second Kind spaces of Harmonic
Dimension $N$ appear in a natural way when we consider the space
$X_{N}\circleddash X_{M}$ based on solutions of Elliptic BVP
(\ref{gorthogonal1})-(\ref{gorthogonal3}).

We give the following definition.

\begin{definition}
\label{DsecondKind} For an integer $M\geq1$ we say that the linear subspace
$X_{M}\subset L_{2}\left(  D\right)  $ is of \textbf{Second Kind} and has
\textbf{Harmonic Dimension} $M,$ and write
\[
\operatorname{hdim}\left(  X_{M}\right)  =M,
\]
if it satisfies all conditions of Definition \ref{Dhdimension} however with an
elliptic operator $P_{2N},$ with $N\geq M$ and all elements $u\in X_{M}$
satisfy $N-M$ boundary conditions
\[
B_{j}u=0\qquad\text{on }\partial D,\ j=1,2,...,N-M.
\]
Here the boundary operators $\left\{  B_{j}\right\}  _{j=1}^{N-M}$ are a
\textbf{normal system} of boundary operators defined on $\partial D,$ by
Definition \ref{Delliptic}, item 4).
\end{definition}

By a technique similar to the already used we may prove the following results
which generalize Theorem \ref{TKolmogorovMultivariate}. We assume that
$K_{p}^{\ast}$ is the set defined by (\ref{KpstarGENERAL}) with a strongly
elliptic \emph{constant coefficients} operator $L_{2p}.$ The space
$\widetilde{X}_{p}$ is defined by (\ref{StildeM}) and the space $\widetilde
{F}_{L}$ by (\ref{FtildeN}).

The following theorem is a generalization of item 1) in Theorem
\ref{TKolmogorovMultivariate}.

\begin{theorem}
\label{TSecondKindMlessp} Let $M<p$ and $L\geq0$ be arbitrary integer. Let
$X_{M}$ be a Second Kind space with Harmonic Dimension $N,$ i.e.
\[
\operatorname{hdim}\left(  X_{M}\right)  =M.
\]
Let $F_{L}$ be an $L-$dimensional subset of $L_{2}\left(  B\right)  .$ Then
\[
\operatorname*{dist}\left(  X_{M}%
{\textstyle\bigoplus}
F_{L},K_{p}^{\ast}\right)  =\infty.
\]

\end{theorem}

The proof of Theorem \ref{TSecondKindMlessp} follows with minor modifications
of Lemma \ref{LtransversalSpaces} (Theorem \ref{TtransversalSpaces}).

It is more non-trivial to consider the case $N=p.$ First we must prove the
following result.

\begin{lemma}
\label{LSecondKindM=p}Let $X_{p}$ be a Second Kind space of Harmonic Dimension
$p$ and $L\geq0$ be an arbitrary integer. Let $F_{L}$ be an $L-$dimensional
subset of $L_{2}\left(  B\right)  .$ Then
\[
\operatorname*{dist}\left(  X_{p}%
{\textstyle\bigoplus}
F_{L},K_{p}^{\ast}\right)  <\infty
\]
implies
\begin{equation}
\widetilde{X}_{p}\subset X_{p}. \label{XptildeConXp}%
\end{equation}
Let the elliptic operator $P_{2M}$ and the boundary operators $\left\{
B_{j}\right\}  _{j=1}^{M-p}$ be associated with $X_{p}$ by Definition
\ref{DsecondKind}. Then (\ref{XptildeConXp}) implies the following
factorizations:
\begin{align*}
P_{2M}  &  =P_{2M-2p}^{\prime}L_{2p}\\
B_{j}  &  =B_{j}^{\prime}L_{2p}\qquad\text{for }j=1,2,...,M-p.
\end{align*}
The operator $P_{2M-2p}^{\prime}$ is uniformly strongly elliptic in $D,$ and
the boundary operators $\left\{  B_{j}^{\prime}\right\}  _{j=1}^{M-p}$ form a
normal system which covers the operator $P_{2M-2p}^{\prime}.$
\end{lemma}

Finally, the following generalization of item 2) in Theorem
\ref{TKolmogorovMultivariate} may be proved. It shows that one needs to take
into account the index of the Elliptic BVP involved.

\begin{theorem}
\label{TKolmogorovSecondKind} Let us consider those spaces $X_{p}$ of Second
Kind with Harmonic Dimension $p$ for which
\[
\operatorname*{dist}\left(  X_{p}%
{\textstyle\bigoplus}
F_{L},K_{p}^{\ast}\right)  <\infty
\]
with associated operators $P_{2M}$ and boundary operators $\left\{
B_{j}\right\}  _{j=1}^{M-p}.$ Following the notations of Lemma
\ref{LSecondKindM=p}, let us denote by $\mathcal{N}$ the following space of
solutions $w\in H^{2M-2p}\left(  D\right)  $ of the Elliptic BVP on the domain
$D$:%
\begin{align*}
P_{2M-2p}^{\prime}w  &  =0\qquad\text{on }D\\
B_{j}^{\prime}w  &  =0,\ \text{on }\partial D,\ \text{for }j=1,2,...,M-p
\end{align*}
Then the following equality holds \
\begin{equation}
\inf_{X_{p},F_{L}}\left\{  \operatorname*{dist}\left(  X_{p}%
{\textstyle\bigoplus}
F_{L},K_{p}^{\ast}\right)  :\text{ }\dim\left(  \mathcal{N}\right)
+L=L_{1}\right\}  =\operatorname*{dist}\left(  \widetilde{X}_{p}%
{\textstyle\bigoplus}
\widetilde{F}_{L_{1}},K_{p}^{\ast}\right)  .\nonumber
\end{equation}

\end{theorem}

From the theory of Elliptic BVP is known that $\dim\left(  \mathcal{N}\right)
<\infty$ (cf. \cite{lions-magenes}, Theorem $5.3,$ chapter $2,$ section
$5.3$). Let us denote by $\left\{  w_{s}\right\}  _{s=1}^{\dim\left(
\mathcal{N}\right)  }$ a basis of the space $\mathcal{N},$ and by $u_{s}$ a
fixed solution to $L_{2p}u_{s}=w_{s}.$ The main point in the proof of Theorem
\ref{TKolmogorovSecondKind} is that arbitrary solution $u$ to equation
$P_{2M}u=0$ may be expressed as
\[
u=%
{\displaystyle\sum_{s=1}^{\dim\left(  \mathcal{N}\right)  }}
\lambda_{s}u_{s}+v
\]
where $v$ is a solution to $L_{2p}v=0.$

\section{Appendix}

\subsection{Proof of Theorem \ref{PstructureMarkovSystem} \label{Sappendix1}}

Follows the proof of Theorem \ref{PstructureMarkovSystem}:%

\proof
Let $\left\{  v_{j}\right\}  _{j=1}^{N}$ be a basis of the space $X_{N},$ i.e.
$X_{N}=\operatorname*{span}\left\{  v_{j}\right\}  _{j=1}^{N}.$ Then by a
Theorem in \cite{Gopengauz} there exists a sequence of algebraic polynomials
$\left\{  P_{j}^{m}\right\}  _{j=1}^{N},$ $m\geq1,$ such that in the norm of
$C^{N-1}\left(  I\right)  $ for $j=1,2,...,N$ holds
\begin{equation}
P_{j}^{m}\longrightarrow v_{j}\qquad\text{for }m\longrightarrow\infty.
\label{Pjmtovj}%
\end{equation}
Since $\dim X_{N}=N$ it is clear that we may choose the polynomials $\left\{
P_{j}^{m}\right\}  _{j=1}^{N}$ also linearly independent. Hence, we obtain the
spaces $X_{N}^{m}=\operatorname*{span}\left\{  P_{j}^{m}\right\}  _{j=1}^{N}$
with $\dim X_{N}^{m}=N.$ Since for the basis we have the limit (\ref{Pjmtovj})
it follows that also for the unit balls holds $U\left(  X_{N}^{m}\right)
\longrightarrow U\left(  X_{N}\right)  ,$ for $m\longrightarrow\infty.$

Let us fix $m\geq1.$ By a theorem of M. B\^{o}cher (cf. \cite{Bostan} and
references therein) we know that the Wronskians of linear independent analytic
functions are never identically zero. Since the Wronskians in (\ref{Wro}) are
given by
\[
W_{k}\left(  t\right)  =W\left(  P_{1}^{m}\left(  t\right)  ,P_{2}^{m}\left(
t\right)  ,...,P_{k}^{m}\left(  t\right)  \right)  \qquad\text{for
}k=1,2,...,N,
\]
they are polynomials, and since they cannot be identically zero it follows
that they have a finite number of zeros. Hence, for fixed $m\geq1$ we may
subdivide the interval $I$ into subintervals $a\leq t_{0}<t_{1}<\cdot
\cdot\cdot<t_{p_{k}}=b,$ and on every subinterval $\left(  t_{j}%
,t_{j+1}\right)  $ all Wronskians have a definite sign. By
\cite{karlinstudden} (chapter $11,$ Theorem $1.1$) it follows that on every
interval $\left(  t_{j},t_{j+1}\right)  $ the polynomials $\left\{  P_{j}%
^{m}\right\}  _{j=1}^{N}$ constitute an $\emph{ECT-}$system. Hence $\left\{
P_{j}^{m}\right\}  _{j=1}^{N}$ is a piecewise $\emph{ECT-}$system. This ends
the proof.%

\endproof

Let us note that if $\left\{  v_{j}\right\}  _{j=1}^{N}$ are polynomials, then
by formula (\ref{roW1})-(\ref{roW2}) the functions $\rho_{j}$ are rational
functions on the whole interval $\left[  a,b\right]  .$ We will provide an
example to illustrate Theorem \ref{PstructureMarkovSystem}.

Let $N=2.$ We consider the $2-$dimensional space $X_{2}=\operatorname*{span}%
\left\{  u_{1},u_{2}\right\}  \subset C^{1}\left(  I\right)  $ where
\begin{align*}
u_{1}\left(  t\right)   &  =\chi\left(  t\right)  t^{2}\\
u_{2}\left(  t\right)   &  =t^{2}%
\end{align*}
and $I=\left(  -1,1\right)  ;$ here $\chi\left(  t\right)  $ is the Heaviside
function equal to $0$ for $t<0$ and to $1$ for $t\geq0.$ Note that the
Wronskian $W\left(  u_{1}\left(  t\right)  ,u_{2}\left(  t\right)  \right)
\equiv0$ on $I.$ However we may approximate the space $X_{2}$ by the linear
$2-$dimensional spaces $X_{2}^{\varepsilon}=\operatorname*{span}\left\{
u_{1}+\varepsilon,u_{2}\right\}  \subset C^{1}\left(  I\right)  $ with
$\varepsilon\longrightarrow0.$ The system $\left\{  u_{1}+\varepsilon
,u_{2}\right\}  $ has a Wronskian $W\left(  u_{1}\left(  t\right)
+\varepsilon,u_{2}\left(  t\right)  \right)  $ equal to $2t\varepsilon.$ For
the last system the functions $\rho_{1}^{\varepsilon}$ and $\rho
_{2}^{\varepsilon}$ satisfying equations (\ref{explicitECT1}%
)-(\ref{explicitECT4}) are given by
\begin{align*}
\rho_{1}^{\varepsilon}\left(  t\right)   &  =u_{1}\left(  t\right)
+\varepsilon\\
\rho_{2}^{\varepsilon}\left(  t\right)   &  =\left\{
\begin{array}
[c]{c}%
\frac{2t}{\varepsilon}\qquad\qquad\text{for }t\leq0\\
\frac{2t\varepsilon}{\left(  t^{2}+\varepsilon\right)  ^{2}}\qquad\text{for
}t\geq0
\end{array}
\right.  .
\end{align*}

\section{Appendix}

\subsection{Proof of Theorem \ref{TExpansionBerezanskii} \label{Sappendix}}

%

\proof
\textbf{(1)} We consider the following auxiliary elliptic \emph{eigenvalue
problem}
\begin{align}
L_{2p}L_{2p}^{\ast}\phi\left(  x\right)   &  =\lambda\phi\left(  x\right)
\qquad\qquad\text{on }B,\label{ffiBVP1}\\
B_{j}\phi\left(  y\right)   &  =S_{j}\phi\left(  y\right)  =0\qquad\text{for
}j=0,1,...,p-1,\text{ for }y\in\partial B. \label{ffiBVP2}%
\end{align}
Since this is the Dirichlet problem for the operator $L_{2p}^{\ast}L_{2p}$ it
is a classical fact that (\ref{ffiBVP1})-(\ref{ffiBVP2}) is a \emph{regular
Elliptic BVP} considered in the Sobolev space $H^{2p}\left(  B\right)  ,$ as
defined in Definition \ref{Delliptic}. Also, it is a classical fact that the
Dirichlet problem is a self-adjoint problem (cf. \cite{lions-magenes}, Remark
$2.4$ in section $2.4$ and Remark $2.6$ in section $2.5$, chapter $2$).

Hence, we may apply the main results about the Spectral theory of regular
self-adjoint Elliptic BVP. We refer to \cite{egorov-shubin0} (section $3$ in
chapter $2,$ p. $122,$ Theorem $2.52$) and to references therein.

By the uniqueness Lemma \ref{LuniquenessDirichlet} the eigenvalue problem
(\ref{ffiBVP1})-(\ref{ffiBVP2}) has only zero solution for $\lambda=0.$ It has
eigenfunctions $\phi_{k}\in H^{2p}\left(  B\right)  $ with eigenvalues
$\lambda_{k}>0$ for $k=1,2,3,...$ for which $\lambda_{k}\longrightarrow\infty$
as $k\longrightarrow\infty.$

\textbf{(2)} Next, in the Sobolev space $H^{2p}\left(  B\right)  ,$ we
consider the problem:
\begin{align}
L_{2p}L_{2p}^{\ast}\varphi\left(  x\right)   &  =\phi_{k}\left(  x\right)
\qquad\qquad\text{on }B\label{fiBVP1}\\
B_{j}\varphi\left(  y\right)   &  =S_{j}\varphi\left(  y\right)
=0\qquad\text{for }j=0,1,...,p-1,\text{ for }y\in\partial B. \label{fiBVP2}%
\end{align}
Obviously, the Elliptic BVP defined by problem (\ref{fiBVP1})-(\ref{fiBVP2})
coincides with the Elliptic BVP defined by (\ref{ffiBVP1})-(\ref{ffiBVP2}) up
to the right-hand sides, and all remarks there hold as well. Hence, problem
(\ref{fiBVP1})-(\ref{fiBVP2}) has\emph{ unique solution} $\varphi_{k}\in
H^{2p}\left(  B\right)  .$ We put
\[
\psi_{k}=L_{2p}^{\ast}\varphi_{k}.
\]
Hence, $L_{2p}\psi_{k}=\phi_{k}.$ We infer that on the boundary $\partial B$
hold the equalities $B_{j}L_{2p}\psi_{k}=B_{j}\phi_{k}$ and $S_{j}L_{2p}%
\psi_{k}=S_{j}\phi_{k};$ since $\phi_{k}$ are solutions to (\ref{ffiBVP1}%
)-(\ref{ffiBVP2}) it follows
\begin{equation}
B_{j}L_{2p}\psi_{k}\left(  y\right)  =S_{j}L_{2p}\psi_{k}\left(  y\right)
=0\qquad\text{for }j=0,1,...,p-1,\text{ for }y\in\partial B. \label{psiBVP}%
\end{equation}

We will prove that $\psi_{k}$ are solutions to problem (\ref{eigen1Multi}%
)-(\ref{eigen2Multi}), they are mutually \textbf{orthogonal,} and they are
also orthogonal to the space \allowbreak\ $\left\{  v\in H^{2p}:L_{2p}%
v=0\right\}  .$

\textbf{(3) }Let us see that
\[
L_{2p}^{\ast}L_{2p}\psi_{k}=\lambda_{k}\psi_{k}.
\]
By the definition of $\psi_{k}$ this is equivalent to
\[
L_{2p}^{\ast}L_{2p}L_{2p}^{\ast}\varphi_{k}=\lambda_{k}L_{2p}^{\ast}%
\varphi_{k};
\]
from $L_{2p}L_{2p}^{\ast}\varphi_{k}=\phi_{k}$ this is equivalent to
\[
L_{2p}^{\ast}\phi_{k}=\lambda_{k}L_{2p}^{\ast}\varphi_{k}%
\]
On the other hand, by the basic properties of $\phi_{k}$ and $\varphi_{k},$ we
have obviously $L_{2p}L_{2p}^{\ast}\phi_{k}=\lambda_{k}L_{2p}L_{2p}^{\ast
}\varphi_{k},$ hence
\[
L_{2p}L_{2p}^{\ast}\left(  \phi_{k}-\lambda_{k}\varphi_{k}\right)  =0.
\]
Note that both $\phi_{k}$ and $\varphi_{k}$ satisfy the same zero Dirichlet
boundary conditions, namely (\ref{ffiBVP2}) and (\ref{fiBVP2}). Hence, by the
uniqueness Lemma \ref{LuniquenessDirichlet} it follows that $\phi_{k}%
-\lambda_{k}\varphi_{k}=0$ which implies $L_{2p}^{\ast}L_{2p}\psi_{k}%
=\lambda_{k}\psi_{k}.$ Thus we see that $\psi_{k}$ is a solution to problem
(\ref{eigen1Multi})-(\ref{eigen2Multi}) and does not satisfy $L_{2p}\psi=0.$

\textbf{(4)} The orthogonality to the subspace $\left\{  v\in H^{2p}%
:L_{2p}v=0\right\}  $ follows easily from the Green formula
(\ref{GreenGeneral}) applied to the operator $L_{2p}^{\ast}L_{2p},$
\begin{align*}
&
{\displaystyle\int_{D}}
\left(  L_{2p}^{\ast}L_{2p}\psi_{k}\cdot v-L_{2p}\psi_{k}\cdot L_{2p}v\right)
dx\\
&  =%
{\displaystyle\sum_{j=0}^{2p-1}}
{\displaystyle\int_{\partial D}}
\left(  S_{j}L_{2p}\psi_{k}\cdot C_{j}v-B_{j}L_{2p}\psi_{k}\cdot
T_{j}v\right)
\end{align*}
in which substitute the zero boundary conditions (\ref{psiBVP}) of $\psi_{k},$
and equality
\[%
{\displaystyle\int_{D}}
L_{2p}^{\ast}L_{2p}\psi_{k}\cdot vdx=\lambda_{k}%
{\displaystyle\int_{D}}
\psi_{k}\cdot vdx.
\]

The orthonormality of the system $\left\{  \psi_{k}\right\}  _{k=1}^{\infty}$
follows now easily by the equality
\[
\lambda_{k}%
{\displaystyle\int}
\psi_{k}\psi_{j}dx=%
{\displaystyle\int}
L_{2p}^{\ast}L_{2p}\psi_{k}\psi_{j}dx=%
{\displaystyle\int}
L_{2p}\psi_{k}L_{2p}\psi_{j}dx=%
{\displaystyle\int}
\phi_{k}\phi_{j}dx
\]
and the orthogonality of the system $\left\{  \phi_{k}\right\}  _{k=1}%
^{\infty}.$

\textbf{(5)} For the completeness of the system $\left\{  \psi_{k}\right\}
_{k=1}^{\infty}$, let us assume that for some $f\in L_{2}\left(  B\right)  $
holds
\begin{equation}%
{\displaystyle\int_{B}}
f\cdot\psi_{k}dx=%
{\displaystyle\int_{B}}
f\cdot\psi_{k}^{\prime}dx=0\qquad\text{for all }k\geq1. \label{fOrthogonal}%
\end{equation}
Then the Green formula (\ref{GreenGeneral}) implies
\begin{align*}
0  &  =\lambda_{k}%
{\displaystyle\int_{B}}
f\cdot\psi_{k}dx=%
{\displaystyle\int_{B}}
f\cdot L_{2p}^{\ast}L_{2p}\psi_{k}dx=%
{\displaystyle\int_{B}}
L_{2p}f\cdot L_{2p}\psi_{k}dx\\
&  =%
{\displaystyle\int_{B}}
L_{2p}f\cdot\phi_{k}dx\qquad\text{for all }k\geq1.
\end{align*}
By the completeness of the system $\left\{  \phi_{k}\right\}  _{k\geq1}$ this
implies that $L_{2p}f=0.$ From the second orthogonality in (\ref{fOrthogonal})
follows that $f\equiv0,$ and this ends the proof of the completeness of the
system $\left\{  \psi_{j}^{\prime}\right\}  _{j=1}^{\infty}\bigcup\left\{
\psi_{j}\right\}  _{j=1}^{\infty}.$%

\endproof

We have used above the following simple result.

\begin{lemma}
\label{LuniquenessDirichlet} The solution to problem (\ref{ffiBVP1}%
)-(\ref{ffiBVP2}) for $\lambda=0$ is unique.
\end{lemma}

%

\proof
From Green's formula (\ref{GreenGeneral}) we obtain
\[%
{\displaystyle\int_{B}}
\left[  L_{2p}\phi\right]  ^{2}dx-%
{\displaystyle\int}
\phi\cdot L_{2p}^{\ast}L_{2p}\phi dx=%
{\displaystyle\sum_{j=1}^{p}}
{\displaystyle\int_{\partial B}}
\left(  S_{j}\phi\cdot C_{j}L_{2p}\phi-B_{j}\phi\cdot T_{j}L_{2p}\phi\right)
d\sigma_{y},
\]
hence $L_{2p}\phi=0.$

Now for arbitrary $v\in H^{2p}\left(  B\right)  $ by the same Green's formula
we obtain
\[%
{\displaystyle\int_{B}}
\left(  L_{2p}\phi\cdot v-\phi\cdot L_{2p}^{\ast}v\right)  dx=%
{\displaystyle\sum_{j=1}^{p}}
{\displaystyle\int_{\partial B}}
\left(  S_{j}\phi\cdot C_{j}v-B_{j}\phi\cdot T_{j}v\right)  d\sigma_{y}=0,
\]
hence
\[%
{\displaystyle\int_{B}}
\phi\cdot L_{2p}^{\ast}vdx=0.
\]
From the local existence theorem for elliptic operators (cf.
\cite{lions-magenes}) it follows that for arbitrary $f\in L_{2}\left(
B\right)  $ we may solve the elliptic equation $L_{2p}^{\ast}v=f$ with $v\in
H^{2p}\left(  B\right)  .$ From the density of $H^{2p}\left(  B\right)  $ in
$L_{2}\left(  B\right)  $ we infer $\phi\equiv0.$

This ends the proof.%

\endproof

\section{Some open problems}

\begin{enumerate}
\item First of all, one has to study basic questions about the sets having
Harmonic Dimension, by considering the sets $X_{M}\bigcap X_{N},$ $X_{M}%
{\textstyle\bigoplus}
X_{N},$ $X_{M}%
{\textstyle\bigotimes}
X_{N},$ and similar, and finding their Harmonic Dimension.

\item One has to check that the maximal generality of the theory in the
present paper will be achieved by considering elliptic pseudo-differential operators.

\item New Jackson type theorems in approximation theory are suggested by the
results proved:\ the simplest way to state them is to consider spaces defined
for example by
\[
\left\{  u:\left\vert L_{2p}u\left(  x\right)  \right\vert \leq1\qquad
\text{for }x\in D\right\}  .
\]
In the case of polyharmonic operator Jackson type results have been proved in
\cite{kounchev1991Hanstholm}. By the proof of Theorem
\ref{TKolmogorovMultivariate} one may expect a Jackson type theorem to be
proved for approximation by solutions of equations $\left\{  u:P_{2N}%
u=0\quad\text{in }D\right\}  ,$ where the operator $P_{2N}$ is of the form
$P_{2N}=P_{2N-2p}^{\prime}L_{2p}.$

\item One has to find a proper discrete version of the present research which
will be essential for the applications to Compressed Sensing, compare the role
of Gelfand's widths in \cite{donoho}, \cite{devore}.

\item Let us recall that one-dimensional Chebyshev systems are important for
the qualitative theory of ODEs, in particular for \emph{Sturmian type}
theorems, cf. e.g. \cite{arnold}, \cite{arnoldProblem}. V.I. Arnold discusses
the importance of the Chebyshev systems in his Toronto lectures, June 1997,
Lecture 3: Topological Problems in Wave Propagation Theory and Topological
Economy Principle in Algebraic Geometry. Fields Institute Communications,
available online at\emph{ http://www.pdmi.ras.ru/\symbol{126}%
arnsem/Arnold/arn-papers.html. } On p. 8 he writes that \textquotedblright
Even the Sturm theory is missing in higher dimensions. This is an interesting
phenomenon. All attempts that I know to extend Sturm theory to higher
dimensions failed. For instance, you can find such an attempt in the
Courant-Hilbert's book, in chapter 6, but it is wrong. The topological
theorems about zeros of linear combinations for higher dimensions, which are
attributed there to Herman, are wrong even for the standard spherical
Laplacian.\textquotedblright\ The attempts to find a proper setting for
multidimensional Chebyshev systems are present in the works of V.I. Arnold in
the context of multivariate Sturm type of theorems, see in particular problem
1996-5 in \cite{arnoldProblem}. In view of this circle of problems, one may
try to apply the present framework and to obtain a proper Ansatz for
multidimensional Sturm type theorems.
\end{enumerate}

\end{document}